\documentclass[openacc]{rsproca_new}
\usepackage{amsmath}
\usepackage{amssymb}
\usepackage{bm}

\newtheorem{lemma}{\bf Lemma}[section]

\def\red{\color{red}}

\begin{document}

\citearticle{M{\"u}ller A. 2021 Review of the exponential and Cayley map on SE(3) as relevant for Lie group integration of the generalized Poisson equation and flexible multibody systems}{20210303}{477, 2021}
\history{{\bf Corrected preprint}\\
Last update: 10 September 2023\\
Corrections are indicated in {\color{red}red}}

%%%% Article title to be placed here
\title{Review of the {Exponential} and Cayley {Map} on SE(3) {as relevant for} Lie Group Integration of the Generalized Poisson Equation and Flexible Multibody Systems}

\author{%%%% Author details
Andreas M\"uller $^{1}$}

%%%%%%%%% Insert author address here
\address{$^{1}$JKU Johannes Kepler University, Linz, Austria}

%%%% Subject entries to be placed here %%%%
\subject{Lie groups, computational mechanics, structural engineering, robotics}

%%%% Keyword entries to be placed here %%%%

%%%% Insert corresponding author and its email address}
\corres{Andreas M\"uller\\
\email{a.mueller@jku.at}}

%%%% Abstract text to be placed here %%%%%%%%%%%%
\begin{abstract}
The exponential and Cayley map on $SE\left( 3\right) $ are the prevailing
coordinate maps used in Lie group integration schemes for rigid body and
flexible body systems. Such geometric integrators are the Munthe-Kaas and
generalized-$\alpha $ schemes, which involve the differential and its
directional derivative of the respective coordinate map. Relevant closed
form expressions, which were reported over the last two decades, are
scattered in the literature, and some are reported without proof. This paper
provides a reference summarizing all relevant closed form relations along
with the relevant proofs. including the right-trivialized differential of
the exponential and Cayley map and their directional derivatives (resembling
the Hessian). The latter gives rise to an implicit generalized-$\alpha $
scheme for rigid/flexible multibody systems in terms of the Cayley map with
improved computational efficiency.
\end{abstract}

\begin{fmtext}
\section{Introduction}

This paper addresses the exponential and Cayley map used to express
solutions of the generalized right and left Poisson equation on $SE\left(
3\right) $ 
\begin{equation}
\dot{\mathbf{C}}=\hat{\mathbf{V}}{^{\mathrm{s}}}\mathbf{C},\ \dot{\mathbf{C}}%
=\mathbf{C}\hat{\mathbf{V}}{^{\mathrm{b}}}  \label{Poisson}
\end{equation}%
{where $\mathbf{C}\left( t\right) \in SE\left( 3\right) $ describes the
motion of a frame in $E^{3}$} (regarded as Euclidean motion), and $\hat{%
\mathbf{V}}{^{\mathrm{s}}}\left( t\right) \in se\left( 3\right) $ is the
velocity in spatial representation, and \mbox{$\hat{\mathbf{V}}{^{\mathrm{b}}}%
\left( t\right) \in se\left( 3\right) $} in body-fixed representation \cite%
{MurrayBook,MUBOScrew1}. This is used to describe the relative motion of
rigid bodies as well as the displacement field of Cosserat continua. \end{fmtext}

\maketitle

The equations (\ref{Poisson})
relate the curve $\mathbf{V}\left( t\right) $ in $se\left( 3\right) $, i.e.
rigid body twist, to the corresponding curve in $SE\left( 3\right) $, i.e.
the rigid body motion. In context of rigid body kinematics and dynamics,
they are thus referred to as the \emph{kinematic reconstruction equations}.
The equations where originally derived by Poisson for pure rotations which
is why (\ref{Poisson}) are referred to as generalized Poisson equations.
Since the equations for pure rotations are also attributed to Darboux \cite%
{Darboux1887}, they are occasionally called Poisson-Darboux equations \cite%
{ConduracheAASAIAA2017}.

\subsection{Lie group integration and coordinate maps}

When seeking solutions of the form $\mathbf{C}\left( t\right) =\exp \hat{%
\mathbf{X}}\left( t\right) \mathbf{C}\left( 0\right) $ in terms of the
canonical coordinates $\mathbf{X}\in {\mathbb{R}}^{6}$ of first kind, then
for small $t $, and with $\hat{\mathbf{X}}\in se\left( 3\right) $ defined in
(\ref{se3}),%
\begin{equation}
\hat{\mathbf{V}}{^{\mathrm{s}}}=\mathrm{dexp}_{\hat{\mathbf{X}}}(\dot{\hat{%
\mathbf{X}}}),\ \ \ \hat{\mathbf{V}}{^{\mathrm{b}}}=\mathrm{dexp}_{-\hat{%
\mathbf{X}}}(\dot{\hat{\mathbf{X}}})  \label{RecExp}
\end{equation}%
where the right-trivialized differential $\mathrm{dexp}_{\hat{\mathbf{X}}%
}:se\left( 3\right) \rightarrow se\left( 3\right) $ is defined by%
\begin{equation}
\left( \mathrm{D}_{\hat{\mathbf{X}}}\exp \right) 
%TCIMACRO{\TeXButton{-0.5ex}{\hspace{-0.5ex}}}%
%BeginExpansion
\hspace{-0.5ex}%
%EndExpansion
(\hat{\mathbf{Y}})=\mathrm{dexp}_{\hat{\mathbf{X}}}(\hat{\mathbf{Y}})\exp (%
\hat{\mathbf{X}})  \label{diff1}
\end{equation}%
with $\mathrm{D}_{\hat{\mathbf{X}}}\exp :se\left( 3\right) \rightarrow
T_{\exp \hat{\mathbf{X}}}SE\left( 3\right) $, so that $\left( \mathrm{D}_{%
\hat{\mathbf{X}}}\exp \right) 
%TCIMACRO{\TeXButton{-0.5ex}{\hspace{-0.5ex}}}%
%BeginExpansion
\hspace{-0.5ex}%
%EndExpansion
(\hat{\mathbf{Y}}):=\frac{d}{dt}\exp (\hat{\mathbf{X}}+t\hat{\mathbf{Y}}%
)|_{t=0}$ is the directional derivative\footnote{%
The directional derivative $\left( D_{\hat{\mathbf{X}}}\exp \right) 
%TCIMACRO{\TeXButton{-0.5ex}{\hspace{-0.5ex}}}%
%BeginExpansion
\hspace{-0.5ex}%
%EndExpansion
(\hat{\mathbf{Y}})$ is also denoted in the literature with $D_{\hat{\mathbf{X%
}}}\exp \cdot \hat{\mathbf{Y}}$.} of $\exp $ at $\hat{\mathbf{X}}$ in
direction of $\hat{\mathbf{Y}}$. This was shown by Magnus \cite{Magnus1954}
for a general Lie group. The second equation in (\ref{RecExp}) simply
follows by changing the sign of $\hat{\mathbf{X}}$ and taking into account (%
\ref{Poisson}). Occasionally, $\mathrm{dexp}_{-\hat{\mathbf{X}}}$, which
satisfy%
\begin{equation}
\left( \mathrm{D}_{\hat{\mathbf{X}}}\exp \right) 
%TCIMACRO{\TeXButton{-0.5ex}{\hspace{-0.5ex}}}%
%BeginExpansion
\hspace{-0.5ex}%
%EndExpansion
(\hat{\mathbf{Y}})=\exp (\hat{\mathbf{X}})\mathrm{dexp}_{-\hat{\mathbf{X}}}(%
\hat{\mathbf{Y}})  \label{diff2}
\end{equation}%
is referred to as the left-trivialized differential, and denoted $\mathrm{%
dexp}_{\hat{\mathbf{X}}}^{L}(\hat{\mathbf{Y}}):=\mathrm{dexp}_{-\hat{\mathbf{%
X}}}(\hat{\mathbf{Y}})$. The equations (\ref{RecExp}) will be called the 
\emph{local reconstruction equations}.

The Cayley map $\mathbf{C}\left( t\right) =\mathrm{cay}\hat{\mathbf{X}}%
\left( t\right) \mathbf{C}\left( 0\right) $ allows expressing the solution
in terms of non-canonical coordinates\footnote{%
For simplicity, the same symbol $\mathbf{X}$ is used for either coordinates
as their meaning is clear from the context.} $\mathbf{X}\in {\mathbb{R}}^{6}$%
. It was shown in \cite{DieleLopezPeluso1998} that a solution parameterized
by the Cayley map satisfies, for small $t$, the following relation holds%
\begin{equation}
\hat{\mathbf{V}}{^{\mathrm{s}}}=\mathrm{dcay}_{\hat{\mathbf{X}}}(\dot{\hat{%
\mathbf{X}}}),\ \ \ \hat{\mathbf{V}}{^{\mathrm{b}}}=\mathrm{dcay}_{-\hat{%
\mathbf{X}}}(\dot{\hat{\mathbf{X}}})  \label{RecCay}
\end{equation}%
where now the right-trivialized differential $\mathrm{dcay}_{\hat{\mathbf{X}}%
}:se\left( 3\right) \rightarrow se\left( 3\right) $ is defined by%
\begin{equation}
\left( \mathrm{D}_{\hat{\mathbf{X}}}\mathrm{cay}\right) \left( \mathbf{Y}%
\right) =\mathrm{dcay}_{\hat{\mathbf{X}}}(\hat{\mathbf{Y}})\mathrm{cay(}\hat{%
\mathbf{X}}).  \label{diffCay}
\end{equation}%
In \cite{OwrenMarthinsen2001}, these results were generalized to arbitrary
coordinate maps $\psi $ on a Lie group. Adopted to the Special Euclidean
group, if $\psi :se\left( 3\right) \rightarrow SE\left( 3\right) $ is a
(local) coordinate map, then a solution $\mathbf{C}\left( t\right) =\psi (%
\hat{\mathbf{X}}\left( t\right) )\mathbf{C}\left( 0\right) $ of (\ref%
{Poisson}) satisfies $\hat{\mathbf{V}}{^{\mathrm{s}}}=\mathrm{d}\psi _{\hat{%
\mathbf{X}}}(\dot{\hat{\mathbf{X}}})$, with $\mathrm{d}\psi _{\hat{\mathbf{X}%
}}(\hat{\mathbf{Y}})$ defined by $\left( \mathrm{D}_{\hat{\mathbf{X}}}\psi
\right) (\hat{\mathbf{Y}})=\mathrm{d}\psi _{\hat{\mathbf{X}}}(\hat{\mathbf{Y}%
})\psi (\hat{\mathbf{X}})$. Such alternative coordinate maps are the $k$%
th-order Cayley map \cite{TsiotrasJunkinsSchaub1997,Tsiotras1998} or the
exponential map in terms of dual quaternions \cite%
{Angeles1998,CohenShohamMMT2017,Selig2010} or the dual number formulation 
\cite{HanBauchau2016} of the vector parameterizations of motion \cite%
{BauchauChoi2003,ArgyrisPoterasu1993}. This paper deals with the exponential
and Cayley map as they are the prevailing coordinate maps used for modeling
and Lie group integration.

Various Lie group integration methods depart from the local reconstruction
equations, as discussed in \cite{IserlesMuntheKaasNrsettZanna2000}. Most of
the time integration schemes assume an explicit ODE system, which amounts to
express (\ref{RecExp}) and (\ref{RecCay}), respectively, as%
\begin{equation}
\dot{\hat{\mathbf{X}}}=\mathrm{dexp}_{\hat{\mathbf{X}}}^{-1}(\hat{\mathbf{V}}%
{^{\mathrm{s}}}),\ \ \dot{\hat{\mathbf{X}}}=\mathrm{dcay}_{\hat{\mathbf{X}}%
}^{-1}(\hat{\mathbf{V}}{^{\mathrm{s}}}),  \label{RecExpCay}
\end{equation}%
and accordingly with negative argument $-\hat{\mathbf{X}}$ for the
body-fixed representation of twists in (\ref{Poisson}). In the original
paper \cite{MuntheKaas-BIT1998,MuntheKaas1999}, Munthe-Kaas used the
exponential map along with (\ref{RecExp}). The Cayley map with (\ref{RecCay}%
) was later used, and the heavy top, with motion evolving on $SO\left(
3\right) $, served as a standard example for Lie group methods \cite%
{EngoMarthinsen1998,EngoMarthinsen2001}. Lie group methods using general
coordinate maps were presented in \cite%
{IserlesMuntheKaasNrsettZanna2000,OwrenMarthinsen2001} and applied to rigid
body systems \cite{CelledoniOwren2003}. An excellent overview can be found
in \cite{IserlesMuntheKaasNrsettZanna2000} and \cite{Owren2018}. While most
integration methods \cite%
{CelledoniOwren2003,TerzeZlatarMueller2012,TerzeMueller2016,TerzeZlatarMueller2015}
adopt explicit integration schemes for solving the vector space ODE (\ref%
{RecExp}) or (\ref{RecCay}), the generalized-$\alpha $ method on vector
spaces (originally derived from the semi-implicit Newmark-schemes \cite%
{Newmark1959}), was reformulated to deal with ODE on Lie groups \cite%
{Makinen2001,Saccon}. These Lie group generalized-$\alpha $ methods are
applied to finite rotations \cite{KryslEndres2005} and to multibody systems
comprising rigid as well as flexible bodies \cite%
{ArnoldBrulsCardona2015,ArnoldHante2017,BrulsCardona2010,BrulsCardonaArnold2012}%
, where rigid body motions as well as the deformation field of flexible
beams are consistently represented as curves in $SE\left( 3\right) $.

\subsection{Lie group modeling of flexible beam kinematics}

The kinematics of flexible beams can be modeled as Euclidean motions. Let $%
\mathbf{C}\left( s\right) ,s\in \left[ 0,L\right] $ describe the
configuration of a beam cross section. The beam kinematics can then be
described by%
\begin{equation}
\mathbf{C}^{\prime }=\hat{\bm{\chi}}{^{\mathrm{s}}}\mathbf{C},\ \mathbf{C}%
^{\prime }=\mathbf{C}\hat{\bm{\chi}}{^{\mathrm{b}}}  \label{deform}
\end{equation}%
where $\hat{\bm{\chi}}:\left[ 0,L\right] \rightarrow se\left( 3\right) $
serve as deformation measure to define the strain field. Borri \& Bottasso 
\cite{BorriBottasso1994a,BorriBottasso1994b} were the first to model the
deformation field in geometrically exact beam models as screw (also called
helicoidal) motions, aiming at a strain invariant formulation. The spatial
(right-invariant) deformation measure $\hat{\bm{\chi}}{^{\mathrm{s}}}$ was
used and referred to as the \emph{base-pole generalized curvature}, and $%
\hat{\bm{\chi}}{^{\mathrm{b}}}$ as the \emph{convected generalized curvature}%
. Therefore, the approach was called base-pole (or fixed-pole) formulation.
Using the base-pole formulation, an invariant conserving integration scheme
was presented in \cite{Borri2001a,Borri2001b}. The material (left-invariant)
deformation measure $\hat{\bm{\chi}}{^{\mathrm{b}}}$ was used by Sonneville 
\cite{SonnevilleCardonaBruls2014,SonnevillePhD} who extended the approach to
geometrically exact beams and shells. The main differences of the right- and
the left-invariant formulations are their invariance properties. Also notice
that in \cite{BorriBottasso1994a,BorriBottasso1994b} the adjoint
representation of $SE\left( 3\right) $ is used. With the helicoidal
approximation, beam deformations are expressed in terms of the exponential
map on $SE\left( 3\right) $, and reconstruction equations analog to (\ref%
{diff1}) and (\ref{diff2}) apply. Since Euclidean motions are screw motions,
this allows exact recovery of beam deformations with locally constant
curvature, i.e. pure bending or helical deformation. In the general
situation, the assumption of a helicoidal deformation field is an
approximation. In this context, the importance of using $SE\left( 3\right) $
instead of $SO\left( 3\right) \times {\mathbb{R}}^{3}$ for rigid body
systems as well as geometrically exact beams and shells must be emphasized,
which was discussed in \cite{MMTCSpace2014,BIT2016} and \cite%
{BottassoBorri1998,Borri2001a,BottassoBorriTrainelli2002,BorriBottassoTrainelli2003,SonnevilleBruls2013,SonnevillePhD}%
, respectively. The modeling of beam kinematics with (\ref{deform}) is
recently used to describe soft robots \cite%
{GraziosoDiGironimoSiciliano2019,OrekhovSimaan2020,RendaArmaniniLebastardCandelierBoyer2020}%
.

\subsection{Contribution of this paper}

Key elements of most Lie group integration schemes are the closed form
expressions for the coordinate maps and their trivialized differentials.
Moreover, the (semi-)implicit generalized-$\alpha $ method, necessitates the
directional derivative of the trivialized differentials in order to
construct the Hessian for the iteration steps involved. Closed form
expressions for relevant relations of the exponential map were already
reported in various publications, and it is difficult to find the
corresponding proofs. While closed form expressions were also reported for
the basic relations of the Cayley map, the differential and its derivative
seem not be published. The motivation of this paper is to provide a
comprehensive reference including all relevant proofs. The reader is
referred to the seminal papers \cite%
{BottassoBorri1998,BorriBottassoTrainelli2003} of Borri \& Bottasso, where
several of the presented relations were already reported using different
notions and approaches. These papers seem not have found their due
recognition.

The paper consists of two main parts. Section \ref{secExp} addresses the
parameterization of motions using the exponential map, while section \ref%
{secCay} deals with the motion representation and parameterization using the
Cayley map. Section \ref{secExp}\ref{secExpSO3} recalls well-known
expressions for rotation parameterization in terms of canonical coordinates
(axis/angle), and its differential along with several relations that are
crucial for structure preserving integration schemes. Section \ref{secExp}%
\ref{secExpSE3} presents closed form relations for the exponential of
Euclidean motions in terms of screw coordinates, and the related expressions
for the trivialized differential and their directional derivatives are
derived. The complete list of closed form relations and the corresponding
proofs is the contribution of section \ref{secExp}. For completeness, the
adjoint representation of $SE\left( 3\right) $ is discussed, which is used
for geometrically exact beam modeling. The Cayley map for representing
rotations and Euclidean motions is recalled in section \ref{secCay}\ref%
{secCaySO3} and \ref{secCay}\ref{secCaySE3}, respectively, and the
trivialized differential and its directional derivative is derived. Finally,
the adjoint representation is considered, and it is shown that the
differentials of the Cayley map of $SE\left( 3\right) $ and of its adjoint
representation are different. The closed form expressions for the
differential of the Cayley map and its derivative are the main contribution
of section \ref{secCay}, which provides the basis for a generalized-$\alpha $
integration method in terms of the Cayley map. The paper closes with a short
conclusion in section \ref{secConclusion}.

\section{Motion Parameterization via the Exponential Map%
%TCIMACRO{\TeXButton{secExp}{\label{secExp}}}%
%BeginExpansion
\label{secExp}%
%EndExpansion
}

The exponential map $\exp :\mathfrak{g}\rightarrow G$ on a $n$-dimensional
Lie group $G$ admits the series expansion \cite{Helgason1979}

\begin{equation}
\exp (\tilde{\mathbf{x}})=\sum_{i=0}^{\infty }\frac{1}{i!}\tilde{\mathbf{x}}%
^{i}.  \label{expSeries}
\end{equation}%
The Lie algebra $\mathfrak{g}$ is assumed isomorphic to ${\mathbb{R}}^{n}$,
and $\tilde{\mathbf{x}}\in \mathfrak{g}$ denotes the matrix associated with
the vector $\mathbf{x}\in {\mathbb{R}}^{n}$. For notational convenience,
throughout the paper, $\tilde{\mathbf{x}}\in so\left( 3\right) $ corresponds
to $\mathbf{x}\in {\mathbb{R}}^{3}$, while $\widehat{\mathbf{X}}\in se\left(
3\right) $ is constructed from $\mathbf{X}\in {\mathbb{R}}^{6}$ \cite%
{MurrayBook}.

The right-trivialized differential $\mathrm{dexp}_{\tilde{\mathbf{x}}}:%
\mathfrak{g}\rightarrow \mathfrak{g}$ of the exponential map on a Lie group
admits the series expansion%
\begin{equation}
\mathrm{dexp}_{\tilde{\mathbf{x}}}\left( \tilde{\mathbf{y}}\right)
=\sum_{i=0}^{\infty }\frac{1}{\left( i+1\right) !}\mathrm{ad}_{\tilde{%
\mathbf{x}}}^{i}\left( \tilde{\mathbf{y}}\right) .  \label{dexpSeries}
\end{equation}%
This was shown by Hausdorff in \cite[pp. 26 \& 36ff]{Hausdorff1906}, and a
proof can be found in \cite{Iserles1984}. The series expansion of the
left-trivialized differential is obtained by inserting $-\mathbf{x}$ in (\ref%
{dexpSeries}), which was also shown in \cite{Hausdorff1906} and a proof is
given in \cite[Theorem 2.14.3.]{Varadarajan1984}. Hausdorff \cite[pp. 27 \&
36ff]{Hausdorff1906} further showed that the inverse of the
right-trivialized differential can be represented by the series%
\begin{equation}
\mathrm{dexp}_{\tilde{\mathbf{x}}}^{-1}\left( \tilde{\mathbf{y}}\right)
=\sum_{i=0}^{\infty }\frac{B_{i}}{i!}\mathrm{ad}_{\tilde{\mathbf{x}}%
}^{i}\left( \tilde{\mathbf{y}}\right)  \label{dexpInvSeries}
\end{equation}%
with the Bernoulli numbers $B_{i}$.

In the following, when representing Lie algebra elements $\tilde{\mathbf{x}}%
\in \mathfrak{g}$ as vectors $\mathbf{x}\in {\mathbb{R}}^{n}$, the
differential and inverse are represented by a matrix denoted $\mathbf{dexp}_{%
\mathbf{x}}:{\mathbb{R}}^{n}\rightarrow {\mathbb{R}}^{n}$ and $\mathbf{dexp}%
_{\mathbf{x}}^{-1}:{\mathbb{R}}^{n}\rightarrow {\mathbb{R}}^{n}$ so that $%
\mathrm{dexp}_{\tilde{\mathbf{x}}}\left( \tilde{\mathbf{y}}\right) =\left( 
\mathbf{dexp}_{\mathbf{x}}\mathbf{y}\right) ^{\sim }$ and $\mathrm{dexp}_{%
\tilde{\mathbf{x}}}^{-1}\left( \tilde{\mathbf{y}}\right) =(\mathbf{dexp}_{%
\mathbf{x}}^{-1}\mathbf{y})^{\sim }$, respectively\footnote{%
In computational mechanics, in particular in context of numerical
integration of multibody systems, the matrix representation of $\mathrm{dexp}
$ is frequently called \emph{tangent operator}.}.

\subsection{Spatial Rotations -- $SO\left( 3\right) $%
%TCIMACRO{\TeXButton{secExpSO3}{\label{secExpSO3}}}%
%BeginExpansion
\label{secExpSO3}%
%EndExpansion
}

\subsubsection{Exponential map --Euler-Rodrigues formula}

The exponential map on $SO\left( 3\right) $ is the analytic form of the
classical result in rigid body kinematics, according to which the rotation
about a given axis can be expressed in closed form, which is attributed to
Euler \cite{Euler1775} and Rodrigues \cite{Rodrigues1840}, see also \cite%
{AltmannBook1986}. For later use introduce the abbreviations 
\begin{align}
\alpha & :=\mathrm{sinc}\varphi =\cos \frac{\varphi }{2}\mathrm{sinc\,}\frac{%
\varphi }{2},\ \ \beta :=\mathrm{sinc}^{2}\mathrm{\,}\frac{\varphi }{2}
\label{abc} \\
\gamma & :=\frac{\alpha }{\beta }=\frac{\cos \frac{\varphi }{2}}{\mathrm{sinc%
}\frac{\varphi }{2}}=\frac{\cos ^{2}\frac{\varphi }{2}}{\mathrm{sinc\,}%
\varphi }=\frac{\mathrm{sinc\,}\varphi }{\mathrm{sinc}^{2}\frac{\varphi }{2}}%
,\ \ \delta :=\frac{1-\mathrm{sinc\,}\varphi }{\varphi ^{2}}=\frac{1-\alpha 
}{\varphi ^{2}}.  \notag
\end{align}%
These terms allow minimizing the number of transcendental function
evaluations.

A skew-symmetric matrix $\tilde{\mathbf{x}}\in so\left( 3\right) $ satisfies
the characteristic equation ${\tilde{\mathbf{x}}^{3}}+\left\Vert \mathbf{x}%
\right\Vert ^{2}\tilde{\mathbf{x}}=\mathbf{0}$. Denote with $\varphi
:=\left\Vert \mathbf{x}\right\Vert $ the rotation angle. The series (\ref%
{expSeries}) gives rise to various closed forms \cite%
{BottemaRoth1979,AltmannBook1986,McCarthyBook1990,MurrayBook}%
\begin{align}
\exp \tilde{\mathbf{x}}& =\mathbf{I}+\tfrac{\sin \left\Vert \mathbf{x}%
\right\Vert \,}{\left\Vert \mathbf{x}\right\Vert }\tilde{\mathbf{x}}+\tfrac{%
1-\cos \left\Vert \mathbf{x}\right\Vert }{\left\Vert \mathbf{x}\right\Vert
^{2}}\,\tilde{\mathbf{x}}^{2}  \label{SO3exp1} \\
& =\mathbf{I}+\alpha \tilde{\mathbf{x}}+\tfrac{1}{2}\beta \tilde{\mathbf{x}}%
^{2}  \label{SO3exp6} \\
\exp \left( \varphi \tilde{\mathbf{n}}\right) & =\mathbf{I}+\sin \varphi \,%
\tilde{\mathbf{n}}+\left( 1-\cos \varphi \right) \,\tilde{\mathbf{n}}^{2} \\
& =\mathbf{I}+\sin \varphi \,\tilde{\mathbf{n}}+%
%TCIMACRO{\TeXButton{red}{\color{red}} }%
%BeginExpansion
\color{red}
%EndExpansion
2\sin ^{2}\frac{\varphi }{2}\,\tilde{\mathbf{n}}^{2}%
%TCIMACRO{\TeXButton{black}{\color{black}}}%
%BeginExpansion
\color{black}%
%EndExpansion
\end{align}%
where $\mathbf{x}:=\varphi \mathbf{n}$, with rotation angle $\varphi $ and
the unit vector $\mathbf{n}$ along the rotation axis. The components of $%
\mathbf{x}\in {\mathbb{R}}^{3}$ serve as canonical coordinates of first
kind. The form (\ref{SO3exp6}) along with the parameters $\alpha ,\beta
,\gamma $ were reported in \cite{ParkChung2005}.

\subsubsection{Differential of the exponential map%
%TCIMACRO{\TeXButton{secDexpSO3}{\label{secDexpSO3}}}%
%BeginExpansion
\label{secDexpSO3}%
%EndExpansion
}

When representing $\tilde{\mathbf{x}}\in so\left( 3\right) $ as vector $%
\mathbf{x}\in {\mathbb{R}}^{3}$, the adjoint operator $\mathrm{ad}_{%
\widetilde{\mathbf{x}}}:so\left( 3\right) \rightarrow so\left( 3\right) $ is
simply ${\tilde{\mathbf{x}}:{\mathbb{R}}^{3}\rightarrow {\mathbb{R}}^{3}}$,
and evaluating (\ref{dexpSeries}) along with ${\tilde{\mathbf{x}}^{3}}%
=-\left\Vert \mathbf{x}\right\Vert ^{2}\tilde{\mathbf{x}}$ yields the
following different matrix forms $\mathbf{dexp}_{\mathbf{x}}:{\mathbb{R}}%
^{3}\rightarrow {\mathbb{R}}^{3}$ of the right-trivialized differential 
\begin{align}
\mathbf{dexp}_{\mathbf{x}}& =\mathbf{I}+\tfrac{1-\cos \left\Vert \mathbf{x}%
\right\Vert }{\left\Vert \mathbf{x}\right\Vert ^{2}}\tilde{\mathbf{x}}+%
\tfrac{\left\Vert \mathbf{x}\right\Vert -\sin \left\Vert \mathbf{x}%
\right\Vert }{\left\Vert \mathbf{x}\right\Vert ^{3}}\tilde{\mathbf{x}}^{2}
\label{dexpSO31} \\
& =\tfrac{1}{\left\Vert \mathbf{x}\right\Vert ^{2}}%
%TCIMACRO{\TeXButton{Big}{\big}}%
%BeginExpansion
\big%
%EndExpansion
[(\mathbf{I}-\exp \tilde{\mathbf{x}})\tilde{\mathbf{x}}+\mathbf{xx}^{T}%
%TCIMACRO{\TeXButton{Big}{\big}}%
%BeginExpansion
\big%
%EndExpansion
]  \label{dexpSO32} \\
& =%
%TCIMACRO{\TeXButton{red}{\color{red}} }%
%BeginExpansion
\color{red}
%EndExpansion
\mathbf{nn}^{T}-\alpha \tilde{\mathbf{n}}^{2}+\tfrac{\beta }{2}\tilde{%
\mathbf{x}}%
%TCIMACRO{\TeXButton{black}{\color{black}} }%
%BeginExpansion
\color{black}
%EndExpansion
\label{dexpSO33} \\
& =\mathbf{I}+\tfrac{1}{2}\mathrm{sinc}^{2}\tfrac{\left\Vert \mathbf{x}%
\right\Vert }{2}\,\tilde{\mathbf{x}}+\tfrac{1}{\left\Vert \mathbf{x}%
\right\Vert ^{2}}\left( 1-\tfrac{1}{2}\mathrm{sinc}\tfrac{\left\Vert \mathbf{%
x}\right\Vert }{2}\cos \tfrac{\left\Vert \mathbf{x}\right\Vert }{2}\right) 
\tilde{\mathbf{x}}^{2}  \label{dexpSO34} \\
& =\mathbf{I}+\tfrac{\beta }{2}\tilde{\mathbf{x}}+\tfrac{1}{\left\Vert 
\mathbf{x}\right\Vert ^{2}}\left( 1-\alpha \right) \tilde{\mathbf{x}}^{2}=%
%TCIMACRO{\TeXButton{red}{\color{red}} }%
%BeginExpansion
\color{red}
%EndExpansion
\mathbf{I}+\tfrac{\beta }{2}\tilde{\mathbf{x}}+\left( 1-\alpha \right) 
\tilde{\mathbf{n}}^{2}%
%TCIMACRO{\TeXButton{black}{\color{black}}}%
%BeginExpansion
\color{black}%
%EndExpansion
.  \label{dexpSO35}
\end{align}%
The form (\ref{dexpSO31}) was reported in \cite{ParkPhD} and later in \cite%
{BulloMurray1995,IbrahimbegovicFreyKozar1995}, but was already derived
before in \cite{BorriMelloAtluri1990}. Expression (\ref{dexpSO35}) was
reported in \cite{ParkChung2005}. These explicit expressions reveal the
singularity at $\left\Vert \mathbf{x}\right\Vert =0$, which reflects the
fact that the 3-dimensional group $SO\left( 3\right) $ is not simply
connected \cite{AltmannBook1986}, and thus does not admit a global
singularity-free parameterization in terms of three canonical coordinates. A
singularity-free parameterization needs at least four (dependent)
parameters. Using unit quaternions (Euler parameters) is a common choice,
which implies using $SU\left( 2\right) $ to represent rotations \cite%
{AltmannBook1986,CND2016}. However, this singularity is removable, i.e. the
limit $\lim_{\left\Vert \mathbf{x}\right\Vert \rightarrow 0}\mathbf{dexp}_{%
\mathbf{x}}=\mathbf{I}$ exists, and the relations (\ref{dexpSO33}) and (\ref%
{dexpSO35}) admit a numerically stable evaluation by replacing the term $%
\tilde{\mathbf{x}}^{2}/\left\Vert \mathbf{x}\right\Vert ^{2}$ with $\tilde{%
\mathbf{n}}^{2}$, and $\mathbf{xx}^{T}/\left\Vert \mathbf{x}\right\Vert ^{2}$
with $\mathbf{nn}^{T}$.

The matrix form of the inverse of the right-trivialized differential, $%
\mathrm{dexp}_{\tilde{\mathbf{x}}}^{-1}:so\left( 3\right) \rightarrow
so\left( 3\right) $, is found from the series (\ref{dexpInvSeries}) as 
\begin{align}
\mathbf{dexp}_{\mathbf{x}}^{-1}& =\mathbf{I}-\tfrac{1}{2}\tilde{\mathbf{x}}+%
\frac{1}{\left\Vert \mathbf{x}\right\Vert ^{2}}\left( 1-\tfrac{\left\Vert 
\mathbf{x}\right\Vert }{2}\cot \tfrac{\left\Vert \mathbf{x}\right\Vert }{2}%
\right) \tilde{\mathbf{x}}^{2}  \label{SO3dexpInv1} \\
& =\mathbf{I}-\tfrac{1}{2}\tilde{\mathbf{x}}+\frac{1}{\left\Vert \mathbf{x}%
\right\Vert ^{2}}\left( 1-\tfrac{\mathrm{sinc}\left\Vert \mathbf{x}%
\right\Vert }{\mathrm{sinc}^{2}\tfrac{\left\Vert \mathbf{x}\right\Vert }{2}}%
\right) \tilde{\mathbf{x}}^{2}  \label{SO3dexpInv2} \\
& =\mathbf{I}-\tfrac{1}{2}\tilde{\mathbf{x}}+\frac{1}{\left\Vert \mathbf{x}%
\right\Vert ^{2}}\left( 1-\tfrac{\cos \frac{\left\Vert \mathbf{x}\right\Vert 
}{2}}{\mathrm{sinc}\frac{\left\Vert \mathbf{x}\right\Vert }{2}}\right) 
\tilde{\mathbf{x}}^{2}=\mathbf{I}-\tfrac{1}{2}\tilde{\mathbf{x}}+\tfrac{1}{%
\left\Vert \mathbf{x}\right\Vert ^{2}}\left( 1-\gamma \right) \tilde{\mathbf{%
x}}^{2}.  \label{SO3dexpInv3}
\end{align}%
The expression (\ref{SO3dexpInv1}) was reported in \cite%
{ParkPhD,BulloMurray1995,IbrahimbegovicFreyKozar1995}, and (\ref{SO3dexpInv3}%
) was used in \cite{ParkChung2005}. In the form (\ref{SO3dexpInv2}) and (\ref%
{SO3dexpInv3}), the inverse can be numerically evaluated without
difficulties.

The matrix of the left-trivialized differential and its differential enjoy
the property $\mathbf{dexp}_{-\mathbf{x}}=\mathbf{dexp}_{\mathbf{x}}^{T}$
and $\mathbf{dexp}_{-\mathbf{x}}^{-1}=\mathbf{dexp}_{\mathbf{x}}^{-T}$,
respectively. The spatial and body-fixed angular velocity are hence related
to the time derivative of the rotation axis times angle via%
\begin{align}
%TCIMACRO{\TeXButton{w}{\bm{\omega }}}%
%BeginExpansion
\bm{\omega }%
%EndExpansion
^{\mathrm{s}} &=\mathbf{dexp}_{\mathbf{x}}\dot{\mathbf{x}},\ \ \dot{\mathbf{x%
}}=\mathbf{dexp}_{\mathbf{x}}^{-1}%
%TCIMACRO{\TeXButton{w}{\bm{\omega }}}%
%BeginExpansion
\bm{\omega }%
%EndExpansion
^{\mathrm{s}} \\
%TCIMACRO{\TeXButton{w}{\bm{\omega }}}%
%BeginExpansion
\bm{\omega }%
%EndExpansion
^{\mathrm{b}} &=\mathbf{dexp}_{\mathbf{x}}^{T}\dot{\mathbf{x}},\ \ \dot{%
\mathbf{x}}=\mathbf{dexp}_{\mathbf{x}}^{-T}%
%TCIMACRO{\TeXButton{w}{\bm{\omega }}}%
%BeginExpansion
\bm{\omega }%
%EndExpansion
^{\mathrm{b}}.
\end{align}

\begin{lemma}
%TCIMACRO{\TeXButton{LemmaSO3dexp}{\label{LemmaSO3dexp}}}%
%BeginExpansion
\label{LemmaSO3dexp}%
%EndExpansion
The exponential map on $SO\left( 3\right) $ and its differential are related
via%
\begin{align}
\exp \tilde{\mathbf{x}}& =\mathbf{dexp}_{\mathbf{x}}^{-T}\mathbf{dexp}_{%
\mathbf{x}}  \label{dexpdexpSO31} \\
& =\mathbf{dexp}_{\mathbf{x}}\mathbf{dexp}_{\mathbf{x}}^{-T}
\label{dexpdexpSO32} \\
& =\mathbf{I}+\tilde{\mathbf{x}}\,\mathbf{dexp}_{\mathbf{x}}
\label{dexpdexpSO33} \\
& =\mathbf{I}+\mathbf{dexp}_{\mathbf{x}}\tilde{\mathbf{x}}\,.
\label{dexpdexpSO34}
\end{align}
\end{lemma}

\begin{proof}
Equating the differential of $\mathbf{R}=\exp \tilde{\mathbf{x}}$ written in
terms of the right- and left-trivialized differential (\ref{diff1}) and (\ref%
{diff2}) yields $\mathrm{dexp}_{\tilde{\mathbf{x}}}\left( \tilde{\mathbf{y}}%
\right) \mathbf{R}=\mathbf{R}\mathrm{dexp}_{-\tilde{\mathbf{x}}}\left( 
\tilde{\mathbf{y}}\right) $, and thus $\mathrm{dexp}_{\tilde{\mathbf{x}}%
}\left( \tilde{\mathbf{y}}\right) =\mathbf{R}\mathrm{dexp}_{-\tilde{\mathbf{x%
}}}\left( \tilde{\mathbf{y}}\right) \mathbf{R}^{T}$, with arbitrary $\mathbf{%
y}\in {\mathbb{R}}^{3}$. Using the matrix form of dexp, this can be written
as $\left( \mathbf{dexp}_{\tilde{\mathbf{x}}}\mathbf{y}\right) ^{\sim
}=\left( \mathbf{R\,dexp}_{-\tilde{\mathbf{x}}}\mathbf{y}\right) ^{^{\sim }}$%
, which yields $\mathbf{R\,dexp}_{-\tilde{\mathbf{x}}}=\mathbf{dexp}_{\tilde{%
\mathbf{x}}}$ and thus (\ref{dexpdexpSO31}). From the explicit form (\ref%
{SO3dexpInv1}) follows that $\tilde{\mathbf{x}}=\mathbf{dexp}_{\tilde{%
\mathbf{x}}}^{-T}-\mathbf{dexp}_{\tilde{\mathbf{x}}}^{-1}$. Multiplication
with $\mathbf{dexp}_{\tilde{\mathbf{x}}}$ from the left yields $\mathbf{I}+\,%
\mathbf{dexp}_{\tilde{\mathbf{x}}}\tilde{\mathbf{x}}=\mathbf{dexp}_{\tilde{%
\mathbf{x}}}\mathbf{dexp}_{\tilde{\mathbf{x}}}^{-T}$ and noting (\ref%
{dexpdexpSO31}) yields (\ref{dexpdexpSO34}). The series expansion (\ref%
{dexpInvSeries}) shows that $\mathbf{dexp}_{\tilde{\mathbf{x}}}\tilde{%
\mathbf{x}}=\tilde{\mathbf{x}}\,\mathbf{dexp}_{\tilde{\mathbf{x}}}$ and
hence (\ref{dexpdexpSO33}). Finally, left-multiplication of $\tilde{\mathbf{x%
}}=\mathbf{dexp}_{\tilde{\mathbf{x}}}^{-T}-\mathbf{dexp}_{\tilde{\mathbf{x}}%
}^{-1}$ with $\mathbf{dexp}_{\tilde{\mathbf{x}}}$, noting (\ref{dexpdexpSO33}%
), shows (\ref{dexpdexpSO32}).
\end{proof}

The relations (\ref{dexpdexpSO31}) - (\ref{dexpdexpSO34}) were reported in 
\cite{BorriMelloAtluri1990} and served as key relations for deriving
conservative integration schemes for flexible systems. Relations (\ref%
{dexpdexpSO31}), (\ref{dexpdexpSO32}) and a slightly different form of (\ref%
{dexpdexpSO33}), (\ref{dexpdexpSO34}) were derived in \cite{BulloMurray1995}
with help of computer algebra software.

\subsection{Euclidean Motions -- $SE\left( 3\right) $%
%TCIMACRO{\TeXButton{secExpSE3}{\label{secExpSE3}}}%
%BeginExpansion
\label{secExpSE3}%
%EndExpansion
}

\subsubsection{Exponential map}

The motion of a rigid body evolves on the Lie group $SE\left( 3\right)
=SO\left( 3\right) \ltimes {\mathbb{R}}^{3}$, which is represented as
subgroup of $GL\left( 4\right) $ with elements%
\begin{equation}
\mathbf{C}=\left( 
\begin{array}{cc}
\mathbf{R} & \mathbf{r} \\ 
\mathbf{0} & 1%
\end{array}%
\right) \in SE\left( 3\right)  \label{C}
\end{equation}%
where $\mathbf{R}\in SO\left( 3\right) $ describes the rotation of a
body-fixed frame, and $\mathbf{r}\in {\mathbb{R}}^{3}$ is the position
vector of the frame origin w.r.t. a reference frame. The Lie algebra $%
se\left( 3\right) $ consists of matrices%
\begin{equation}
\hat{\mathbf{X}}=\left( 
\begin{array}{cc}
\tilde{\mathbf{x}} & \mathbf{y} \\ 
\mathbf{0} & 0%
\end{array}%
\right) \in se\left( 3\right)  \label{se3}
\end{equation}%
with $\tilde{\mathbf{x}}\in so\left( 3\right) $ and $\mathbf{y}\in {\mathbb{R%
}}^{3}$. There is an obvious correspondence of $\hat{\mathbf{X}}\in se\left(
3\right) $ and $\mathbf{X}=\left( \mathbf{x},\mathbf{y}\right) \in {\mathbb{R%
}}^{6}$, which renders $se\left( 3\right) $ isomorphic to ${\mathbb{R}}^{6}$.

Chasles theorem \cite{Chasles1830,Chasles1843} asserts that any \emph{finite
rigid body displacement} $\mathbf{C}\in SE\left( 3\right) $ can be achieved
by a screw motion about a fixed screw axis, i.e. there is a constant $%
\mathbf{X}$ so that the screw motion is $\exp (t\hat{\mathbf{X}}),t\in
\lbrack 0,1]$ and $\mathbf{C}=\exp \hat{\mathbf{X}}$. According to a theorem
by Mozzi \cite{Mozzi1763} and Cauchy \cite{Cauchy1827}, for any \emph{given
motion} in time $\mathbf{C}\left( t\right) $, there is an instantaneous
screw axis, i.e. there is a $\mathbf{X}\left( t\right) $ so that $\mathbf{C}%
\left( t\right) =\exp \hat{\mathbf{X}}\left( t\right) \mathbf{C}\left(
0\right) $. The components of $\mathbf{X}$ are the \emph{screw coordinates} 
\cite{SeligBook,MurrayBook}. Solving the kinematic reconstruction equations
thus amounts to determine $\mathbf{X}\left( t\right) $ for given twist.

Evaluating the series (\ref{expSeries}) with matrices (\ref{se3}) yields%
\begin{align}
\exp (\hat{\mathbf{X}})& =\left( 
\begin{array}{cc}
\sum_{i=0}^{\infty }\frac{1}{i!}\tilde{\mathbf{x}}^{i} & \ \
\sum_{i=0}^{\infty }\frac{1}{\left( i+1\right) !}\tilde{\mathbf{x}}^{i}%
\mathbf{y} \\ 
\mathbf{0} & 1%
\end{array}%
\right)  \notag \\
& =\left( 
\begin{array}{cc}
\exp \tilde{\mathbf{x}} & \ \ \mathbf{dexp}_{\mathbf{x}}\mathbf{y} \\ 
\mathbf{0} & 1%
\end{array}%
\right) .  \label{expSE3}
\end{align}%
The dexp mapping on $SO\left( 3\right) $ occurs since $SE\left( 3\right) $
is the semidirect product with ${\mathbb{R}}^{3}$ as normal subgroup, and $%
se\left( 3\right) $ is the semidirect sum with ${\mathbb{R}}^{3}$ as ideal
(i.e., translations do not affect rotations). Inserting (\ref{dexpSO32})
into (\ref{expSE3}) yields%
\begin{align}
\exp (\hat{\mathbf{X}})& =\left( 
\begin{array}{cc}
\mathbf{R} & \ \ \tfrac{1}{\left\Vert \mathbf{x}\right\Vert ^{2}}(\mathbf{I}-%
\mathbf{R})\tilde{\mathbf{x}}\mathbf{y}+h\mathbf{y} \\ 
\mathbf{0} & 1%
\end{array}%
\right) ,\ \mathrm{for\ }\mathbf{x}\neq \mathbf{0}\   \notag \\
& =\left( 
\begin{array}{cc}
\mathbf{I} & \ \mathbf{y} \\ 
\mathbf{0} & 1%
\end{array}%
\right) ,\ \mathrm{for\ }\mathbf{x}=\mathbf{0}  \label{SE3exp1}
\end{align}%
with rotation matrix $\mathbf{R}=\exp \tilde{\mathbf{x}}$. The term $h:=%
\mathbf{x}^{T}\mathbf{y}/\left\Vert \mathbf{x}\right\Vert ^{2}$ is the
instantaneous pitch of the screw motion. For pure rotation $h=0$, and a pure
translations corresponds to $\mathbf{X}=\left( \mathbf{0},\mathbf{y}\right) $%
, i.e. $\mathbf{x}=\mathbf{0}$, for which $h=\infty $. The translation is
then described by the vector $\mathbf{y}=\varphi \mathbf{n}$, where $\mathbf{%
n}\in {\mathbb{R}}^{3}$ is the unit vector in direction of the translation
and $\varphi $ is the amount of translation.

Using (\ref{SO3exp6}), the relation (\ref{SE3exp1}) for $\mathbf{x}\neq 
\mathbf{0}$ can be written as%
\begin{equation}
\exp (\hat{\mathbf{X}})=\left( 
\begin{array}{cc}
\mathbf{R} & \ \ \ 
%TCIMACRO{\TeXButton{red}{\color{red}} }%
%BeginExpansion
\color{red}
%EndExpansion
(\mathbf{I}+\tfrac{\beta }{2}\tilde{\mathbf{x}}+\left( 1-\alpha \right) 
\tilde{\mathbf{n}}^{2})\mathbf{y}%
%TCIMACRO{\TeXButton{black}{\color{black}} }%
%BeginExpansion
\color{black}
%EndExpansion
\\ 
\mathbf{0} & 1%
\end{array}%
\right) ,\ \mathrm{for\ }\mathbf{x}\neq \mathbf{0}\ 
\end{equation}%
with $\mathbf{n}=\mathbf{x/}\left\Vert \mathbf{x}\right\Vert $, which admit
robust numerical evaluation. The exponential can be expressed in terms of
geometric attributes (direction, position, magnitude, pitch) of the
instantaneous screw motion\cite{MUBOScrew1,CND2016,MurrayBook}, which is
relevant for modeling mechanisms and multibody systems in terms of joint
variables and joint screw coordinates.

\subsubsection{Differential of the exponential map}

The velocity of a rigid body in spatial representation, called \emph{spatial
twist}, is defined by $\hat{\mathbf{V}}{^{\mathrm{s}}}=\dot{\mathbf{C}}%
\mathbf{C}^{-1}\in se\left( 3\right) $, where the twist vector $\mathbf{V}^{%
\mathrm{s}}=(%
%TCIMACRO{\TeXButton{w}{\bm{\omega }}}%
%BeginExpansion
\bm{\omega }%
%EndExpansion
^{\mathrm{s}},\mathbf{v}^{\mathrm{s}})\in {\mathbb{R}}^{6}$ comprises the
spatial angular velocity defined by $\tilde{%
%TCIMACRO{\TeXButton{w}{\bm{\omega }}}%
%BeginExpansion
\bm{\omega }%
%EndExpansion
}_{i}^{\mathrm{s}}=\dot{\mathbf{R}}_{i}\mathbf{R}_{i}^{T}$ and spatial
velocity $\mathbf{v}_{i}^{\mathrm{s}}=\dot{\mathbf{r}}_{i}+\mathbf{r}%
_{i}\times 
%TCIMACRO{\TeXButton{w}{\bm{\omega }}}%
%BeginExpansion
\bm{\omega }%
%EndExpansion
_{i}^{\mathrm{s}}$ \cite{MUBOScrew1,MUBOScrew2,MurrayBook}. The velocity in
body-fixed representation, called \emph{body-fixed twist}, is defined by $%
\hat{\mathbf{V}}{^{\mathrm{b}}=}\mathbf{C}^{-1}\dot{\mathbf{C}}\in se\left(
3\right) $, and the twist vector $\mathbf{V}^{\mathrm{b}}=(%
%TCIMACRO{\TeXButton{w}{\bm{\omega }}}%
%BeginExpansion
\bm{\omega }%
%EndExpansion
^{\mathrm{b}},\mathbf{v}^{\mathrm{b}})\in {\mathbb{R}}^{6}$ consists of the
body-fixed angular velocity $\tilde{%
%TCIMACRO{\TeXButton{w}{\bm{\omega }}}%
%BeginExpansion
\bm{\omega }%
%EndExpansion
}_{i}^{\mathrm{b}}=\mathbf{R}_{i}^{T}\dot{\mathbf{R}}_{i}$ and the velocity $%
\mathbf{v}_{i}^{\mathrm{b}}=\mathbf{R}_{i}^{T}\dot{\mathbf{r}}_{i}$.

Expressing the twist in terms of canonical (screw) coordinates $\mathbf{X}$,
according to (\ref{RecExp}), involves the differential $\mathrm{dexp}_{\hat{%
\mathbf{X}}}:se\left( 3\right) \rightarrow se\left( 3\right) $. When
representing twists and screws as vectors, the spatial and body-fixed twist
is determined as ${\mathbf{V}}{^{\mathrm{s}}}=\mathbf{dexp}_{\mathbf{X}}\dot{%
\mathbf{X}}$ and ${\mathbf{V}}{^{\mathrm{b}}}=\mathbf{dexp}_{-\mathbf{X}}%
\dot{\mathbf{X}}$, respectively, where $\mathbf{dexp}_{\mathbf{X}}:{\mathbb{R%
}}^{6}\rightarrow {\mathbb{R}}^{6}$ is the matrix form of the
right-trivialized differential so that $\mathrm{dexp}_{\hat{\mathbf{X}}}(%
\hat{\mathbf{Y}})=\left( \mathbf{dexp}_{\mathbf{X}}\mathbf{Y}\right)
^{\wedge }$. The differential can be determined with the general relation (%
\ref{dexpSeries}) in terms of $\mathrm{ad}_{\hat{\mathbf{X}}}:se\left(
3\right) \rightarrow se\left( 3\right) $. The adjoint operator on $se\left(
3\right) $ defines the Lie bracket as matrix commutator $\mathrm{ad}_{\hat{%
\mathbf{X}}_{1}}(\hat{\mathbf{X}}_{2})=[\hat{\mathbf{X}}_{1},\hat{\mathbf{X}}%
_{2}]=\hat{\mathbf{X}}_{1}\hat{\mathbf{X}}_{2}-\hat{\mathbf{X}}_{2}\hat{%
\mathbf{X}}_{1}$. In vector representation $\mathbf{X}=\left( \mathbf{x},%
\mathbf{y}\right) \in {\mathbb{R}}^{6}$, the matrix form of the adjoint
operator is \cite{SeligBook,MurrayBook}%
\begin{equation}
\mathbf{ad}_{\mathbf{X}}=\left( 
\begin{array}{cc}
\tilde{\mathbf{x}} & \ \mathbf{0} \\ 
\tilde{\mathbf{y}} & \tilde{\mathbf{x}}%
\end{array}%
\right)  \label{adSE3}
\end{equation}%
so that $\hat{\mathbf{Z}}=\mathrm{ad}_{\hat{\mathbf{X}}_{1}}(\hat{\mathbf{X}}%
_{2})$ is equivalent to $\mathbf{Z}=\mathbf{ad}_{\mathbf{X}_{1}}\mathbf{X}%
_{2}$. The result $\mathbf{Z}$ of this operation is also called the screw
product of $\mathbf{X}_{1}$ and $\mathbf{X}_{2}$ \cite%
{McCarthyBook1990,SeligBook}.

\begin{lemma}
%TCIMACRO{\TeXButton{LemmaSE3dexp}{\label{LemmaSE3dexp}}}%
%BeginExpansion
\label{LemmaSE3dexp}%
%EndExpansion
The right-trivialized differential of the exponential map on $SE\left(
3\right) $ possesses the explicit matrix form, with canonical coordinates $%
\mathbf{X}=\left( \mathbf{x},\mathbf{y}\right) $,%
\begin{equation}
\mathbf{dexp}_{\mathbf{X}}=\left( 
\begin{array}{cc}
\mathbf{dexp}_{\mathbf{x}} & \ \ \mathbf{0} \\ 
\left( \mathrm{D}_{\mathbf{x}}\mathbf{dexp}\right) 
%TCIMACRO{\TeXButton{-0.5ex}{\hspace{-0.5ex}}}%
%BeginExpansion
\hspace{-0.5ex}%
%EndExpansion
\left( \mathbf{y}\right) & \mathbf{dexp}_{\mathbf{x}}%
\end{array}%
\right)  \label{dexpSE31}
\end{equation}%
where $\mathbf{dexp}_{\mathbf{x}}$ is the matrix of the right-trivialized
differential on $so\left( 3\right) $ and its directional derivative is given
explicitly as%
\begin{align}
\left( \mathrm{D}_{\mathbf{x}}\mathbf{dexp}\right) 
%TCIMACRO{\TeXButton{-0.5ex}{\hspace{-0.5ex}}}%
%BeginExpansion
\hspace{-0.5ex}%
%EndExpansion
\left( \mathbf{y}\right) =& \tfrac{1}{\left\Vert \mathbf{x}\right\Vert ^{2}}%
\left( 1-\cos \left\Vert \mathbf{x}\right\Vert \right) \tilde{\mathbf{y}}+%
\tfrac{1}{\left\Vert \mathbf{x}\right\Vert ^{3}}\left( \left\Vert \mathbf{x}%
\right\Vert -\sin \left\Vert \mathbf{x}\right\Vert \right) \left( \tilde{%
\mathbf{x}}\tilde{\mathbf{y}}+\tilde{\mathbf{y}}\tilde{\mathbf{x}}\right)
\label{diffDexpSO34} \\
& +\tfrac{\mathbf{x}^{T}\mathbf{y}}{\left\Vert \mathbf{x}\right\Vert ^{4}}%
\left( \left\Vert \mathbf{x}\right\Vert \sin \left\Vert \mathbf{x}%
\right\Vert +2\cos \left\Vert \mathbf{x}\right\Vert -2\right) \tilde{\mathbf{%
x}}+\tfrac{\mathbf{x}^{T}\mathbf{y}}{\left\Vert \mathbf{x}\right\Vert ^{5}}%
\left( 3\sin \left\Vert \mathbf{x}\right\Vert -\left\Vert \mathbf{x}%
\right\Vert \cos \left\Vert \mathbf{x}\right\Vert -2\left\Vert \mathbf{x}%
\right\Vert \right) \tilde{\mathbf{x}}^{2}  \notag \\
=& \tfrac{1}{2}\mathrm{sinc}^{2}\tfrac{\left\Vert \mathbf{x}\right\Vert }{2}%
\,\tilde{\mathbf{y}}+\tfrac{1}{\left\Vert \mathbf{x}\right\Vert ^{2}}\left(
1-\mathrm{\cos }\tfrac{\left\Vert \mathbf{x}\right\Vert }{2}\mathrm{sinc}%
\tfrac{\left\Vert \mathbf{x}\right\Vert }{2}\right) \left( \tilde{\mathbf{x}}%
\tilde{\mathbf{y}}+\tilde{\mathbf{y}}\tilde{\mathbf{x}}\right)
\label{diffDexpSO31} \\
& +\tfrac{\mathbf{x}^{T}\mathbf{y}}{\left\Vert \mathbf{x}\right\Vert ^{2}}%
\left( \mathrm{\cos }\tfrac{\left\Vert \mathbf{x}\right\Vert }{2}\mathrm{sinc%
}\tfrac{\left\Vert \mathbf{x}\right\Vert }{2}-\mathrm{sinc}^{2}\tfrac{%
\left\Vert \mathbf{x}\right\Vert }{2}\right) \tilde{\mathbf{x}}  \notag \\
& +\tfrac{\mathbf{x}^{T}\mathbf{y}}{\left\Vert \mathbf{x}\right\Vert ^{2}}%
\left( \frac{1}{2}\mathrm{sinc}^{2}\tfrac{\left\Vert \mathbf{x}\right\Vert }{%
2}-\tfrac{3}{\left\Vert \mathbf{x}\right\Vert ^{2}}(1-\cos \tfrac{\left\Vert 
\mathbf{x}\right\Vert }{2}\mathrm{sinc}\tfrac{\left\Vert \mathbf{x}%
\right\Vert }{2}\right) \tilde{\mathbf{x}}^{2}  \notag \\
=& \tfrac{1}{2}\mathrm{sinc}^{2}\tfrac{\left\Vert \mathbf{x}\right\Vert }{2}%
\tilde{\mathbf{y}}+\tfrac{1}{\left\Vert \mathbf{x}\right\Vert ^{2}}\left( 1-%
\mathrm{sinc}\left\Vert \mathbf{x}\right\Vert \right) \left( \tilde{\mathbf{x%
}}\tilde{\mathbf{y}}+\tilde{\mathbf{y}}\tilde{\mathbf{x}}\right)
\label{diffDexpSO32} \\
& +\tfrac{\mathbf{x}^{T}\mathbf{y}}{\left\Vert \mathbf{x}\right\Vert ^{2}}%
\left( \mathrm{sinc}\left\Vert \mathbf{x}\right\Vert -\mathrm{sinc}^{2}%
\tfrac{\left\Vert \mathbf{x}\right\Vert }{2}\right) \tilde{\mathbf{x}}+%
\tfrac{\mathbf{x}^{T}\mathbf{y}}{\left\Vert \mathbf{x}\right\Vert ^{2}}%
\left( \tfrac{1}{2}\mathrm{sinc}^{2}\tfrac{\left\Vert \mathbf{x}\right\Vert 
}{2}-\tfrac{3}{\left\Vert \mathbf{x}\right\Vert ^{2}}(1-\mathrm{sinc}%
\left\Vert \mathbf{x}\right\Vert )\right) \tilde{\mathbf{x}}^{2}  \notag \\
=& \tfrac{\beta }{2}\tilde{\mathbf{y}}+\delta \left( \tilde{\mathbf{x}}%
\tilde{\mathbf{y}}+\tilde{\mathbf{y}}\tilde{\mathbf{x}}\right) +\tfrac{%
\mathbf{x}^{T}\mathbf{y}}{\left\Vert \mathbf{x}\right\Vert ^{2}}\left(
\alpha -\beta \right) \tilde{\mathbf{x}}+\tfrac{\mathbf{x}^{T}\mathbf{y}}{%
\left\Vert \mathbf{x}\right\Vert ^{2}}\left( \tfrac{\beta }{2}-3\delta
\right) \tilde{\mathbf{x}}^{2}.  \label{diffDexpSO33}
\end{align}
\end{lemma}

\begin{proof}
Evaluating the series (\ref{dexpSeries}) with (\ref{adSE3}) involves the
powers%
\begin{equation}
\mathbf{ad}_{\mathbf{X}}^{i}=\left( 
\begin{array}{cc}
\tilde{\mathbf{x}}^{i} & \ \mathbf{0} \\ 
\mathbf{P}_{i} & \tilde{\mathbf{x}}^{i}%
\end{array}%
\right) \ \ \ \mathrm{with}\ \ \mathbf{P}_{i}\left( \tilde{\mathbf{x}},%
\tilde{\mathbf{y}}\right) =\sum\limits_{j=0}^{i-1}\tilde{\mathbf{x}}^{j}%
\tilde{\mathbf{y}}\tilde{\mathbf{x}}^{i-j-1},i\geq 1  \label{adi}
\end{equation}%
and $\mathbf{P}_{0}=\mathbf{0}$. The matrix $\mathbf{P}_{i}$ can be written
as directional derivative of the $i$th power $\tilde{\mathbf{x}}^{i}$,
considered as function of $\tilde{\mathbf{x}}$, in the direction of $\tilde{%
\mathbf{y}}$ 
\begin{align}
\mathbf{P}_{i}\left( \tilde{\mathbf{x}},\tilde{\mathbf{y}}\right) =& (%
\mathrm{D}_{\tilde{\mathbf{x}}}\tilde{\mathbf{x}}^{i})%
%TCIMACRO{\TeXButton{-0.5ex}{\hspace{-0.5ex}}}%
%BeginExpansion
\hspace{-0.5ex}%
%EndExpansion
\left( \tilde{\mathbf{y}}\right) =\frac{d}{dt}\left. \left( \tilde{\mathbf{x}%
}+t\tilde{\mathbf{y}}\right) ^{i}\right\vert _{t=0}  \notag \\
=& \left( \tilde{\mathbf{x}}+t\tilde{\mathbf{y}}\right) ^{i-1}\tilde{\mathbf{%
y}}+\left( \tilde{\mathbf{x}}+t\tilde{\mathbf{y}}\right) ^{i-2}\tilde{%
\mathbf{y}}\left( \tilde{\mathbf{x}}+t\tilde{\mathbf{y}}\right) +\ldots
+\left. \left( \tilde{\mathbf{x}}+t\tilde{\mathbf{y}}\right) \tilde{\mathbf{y%
}}^{i-1}\right\vert _{t=0}.  \label{P}
\end{align}%
The series expansion (\ref{dexpSeries}) for the differential thus leads to%
\begin{equation}
\mathbf{dexp}_{\mathbf{X}}=\sum_{i=0}^{\infty }\frac{1}{\left( i+1\right) !}%
\mathbf{ad}_{\mathbf{X}}^{i}=\left( 
\begin{array}{cc}
\mathbf{dexp}_{\mathbf{x}} & \ \ \mathbf{0} \\ 
\sum_{i=0}^{\infty }\frac{1}{\left( i+1\right) !}\left( \mathrm{D}_{\tilde{%
\mathbf{x}}}\tilde{\mathbf{x}}^{i}\right) 
%TCIMACRO{\TeXButton{-0.5ex}{\hspace{-0.5ex}}}%
%BeginExpansion
\hspace{-0.5ex}%
%EndExpansion
\left( \tilde{\mathbf{y}}\right) & \mathbf{dexp}_{\mathbf{x}}%
\end{array}%
\right) =\left( 
\begin{array}{cc}
\mathbf{dexp}_{\mathbf{x}} & \ \ \mathbf{0} \\ 
\left( \mathrm{D}_{\mathbf{x}}\mathbf{dexp}\right) 
%TCIMACRO{\TeXButton{-0.5ex}{\hspace{-0.5ex}}}%
%BeginExpansion
\hspace{-0.5ex}%
%EndExpansion
\left( \mathbf{y}\right) & \mathbf{dexp}_{\mathbf{x}}%
\end{array}%
\right)
\end{equation}%
with the directional derivative of the right-trivialized differential on $%
SO\left( 3\right) $ 
\begin{equation}
\left( \mathrm{D}_{\mathbf{x}}\mathbf{dexp}\right) 
%TCIMACRO{\TeXButton{-0.5ex}{\hspace{-0.5ex}}}%
%BeginExpansion
\hspace{-0.5ex}%
%EndExpansion
\left( \mathbf{y}\right) =\frac{d}{dt}\left. \mathbf{dexp}\left( \tilde{%
\mathbf{x}}+t\tilde{\mathbf{y}}\right) \right\vert _{t=0}.
\end{equation}%
This is evaluated using the following expressions%
\begin{equation}
\left( \mathrm{D}_{\mathbf{x}}\left\Vert \mathbf{x}\right\Vert \right)
\left( \mathbf{y}\right) =(\mathbf{x}^{T}\mathbf{y})/\left\Vert \mathbf{x}%
\right\Vert ,\left( \mathrm{D}_{\mathbf{x}}\tilde{\mathbf{x}}\right) \left( 
\tilde{\mathbf{y}}\right) =\tilde{\mathbf{y}},(\mathrm{D}_{\mathbf{x}}\tilde{%
\mathbf{x}}^{2})\left( \mathbf{y}\right) =\tilde{\mathbf{x}}\tilde{\mathbf{y}%
}+\tilde{\mathbf{y}}\tilde{\mathbf{x}}.  \label{DD}
\end{equation}%
Starting from (\ref{dexpSO31}) yields (\ref{diffDexpSO34}), using (\ref%
{dexpSO34}) yields (\ref{diffDexpSO31}), while using (\ref{dexpSO35}) yields
(\ref{diffDexpSO32}), which can both be expressed as (\ref{diffDexpSO33}).
The closed form (\ref{dexpSE31}) for the differential on $SE\left( 3\right) $
can now be written in terms of these expressions.
\end{proof}

The expression (\ref{diffDexpSO33}) was presented in \cite{ParkChung2005}
without proof in terms of the parameters $\alpha $ and $\beta $. Prior to
this, it was reported in \cite{BorriMelloAtluri1990} and \cite%
{IbrahimbegovicFreyKozar1995} in almost the same form. The singularity at $%
\left\Vert \mathbf{x}\right\Vert =0$ is inherited from the dexp map on $%
SO\left( 3\right) $. The term $\tilde{\mathbf{x}}(\mathbf{x}^{T}\mathbf{y)}%
/\left\Vert \mathbf{x}\right\Vert ^{2}$ can be computed robustly as $\tilde{%
\mathbf{n}}(\mathbf{n}^{T}\mathbf{y})$. The following expression was also
presented in \cite{ParkChung2005}.

\begin{lemma}
The inverse of the right-trivialized differential of the exponential map on $%
SE\left( 3\right) $ possess the matrix form%
\begin{equation}
\mathbf{dexp}_{\mathbf{X}}^{-1}=\left( 
\begin{array}{cc}
\mathbf{dexp}_{\mathbf{x}}^{-1} & \ \ \mathbf{0} \\ 
(\mathrm{D}_{\mathbf{x}}\mathbf{dexp}^{-1})%
%TCIMACRO{\TeXButton{-0.6ex}{\hspace{-0.6ex}}}%
%BeginExpansion
\hspace{-0.6ex}%
%EndExpansion
\left( \mathbf{y}\right) & \mathbf{dexp}_{\mathbf{x}}^{-1}%
\end{array}%
\right)  \label{DexpInvSE3}
\end{equation}%
with canonical coordinates $\mathbf{X}=\left( \mathbf{x},\mathbf{y}\right) $%
, the matrix in (\ref{SO3dexpInv1})-(\ref{SO3dexpInv3}), and 
\begin{equation}
(\mathrm{D}_{\mathbf{x}}\mathbf{dexp}^{-1})%
%TCIMACRO{\TeXButton{-0.5ex}{\hspace{-0.5ex}}}%
%BeginExpansion
\hspace{-0.5ex}%
%EndExpansion
\left( \mathbf{y}\right) =-\tfrac{1}{2}\tilde{\mathbf{y}}+\frac{1}{%
\left\Vert \mathbf{x}\right\Vert ^{2}}\left( 1-\gamma \right) \left( \tilde{%
\mathbf{x}}\tilde{\mathbf{y}}+\tilde{\mathbf{y}}\tilde{\mathbf{x}}\right) +%
\tfrac{\mathbf{x}^{T}\mathbf{y}}{\left\Vert \mathbf{x}\right\Vert ^{4}}%
\left( \tfrac{1}{\beta }+\gamma -2\right) \tilde{\mathbf{x}}^{2}.
\label{diffDexpInvSO3}
\end{equation}
\end{lemma}

\begin{proof}
Invoking the series expansion (\ref{dexpInvSeries}), along with the power of 
$\hat{\mathbf{X}}$ in (\ref{adi}), yields%
\begin{eqnarray}
\mathbf{dexp}_{\mathbf{X}}^{-1} &=&\sum_{i=0}^{\infty }\frac{B_{i}}{i!}%
\mathbf{ad}_{\mathbf{X}}^{i}=\left( 
\begin{array}{cc}
\mathbf{dexp}_{\mathbf{x}}^{-1} & \ \ \mathbf{0} \\ 
\sum_{i=0}^{\infty }\frac{B_{i}}{i!}(\mathrm{D}_{\tilde{\mathbf{x}}}\tilde{%
\mathbf{x}}^{i})%
%TCIMACRO{\TeXButton{-0.7ex}{\hspace{-0.7ex}}}%
%BeginExpansion
\hspace{-0.7ex}%
%EndExpansion
\left( \tilde{\mathbf{y}}\right) & \mathbf{dexp}_{\mathbf{x}}^{-1}%
\end{array}%
\right) \\
&=&\left( 
\begin{array}{cc}
\mathbf{dexp}_{\mathbf{x}}^{-1} & \ \ \mathbf{0} \\ 
(\mathrm{D}_{\mathbf{x}}\mathbf{dexp}^{-1})%
%TCIMACRO{\TeXButton{-0.6ex}{\hspace{-0.6ex}}}%
%BeginExpansion
\hspace{-0.6ex}%
%EndExpansion
\left( \mathbf{y}\right) & \mathbf{dexp}_{\mathbf{x}}^{-1}%
\end{array}%
\right)
\end{eqnarray}%
where the relation (\ref{P}) along with (\ref{dexpInvSeries}) has been used
to identify $\mathrm{D}_{\mathbf{x}}\mathbf{dexp}^{-1}$ as the directional
derivative of $\mathbf{dexp}^{-1}$ on $SO\left( 3\right) $. The statement
follows with the directional derivatives of (\ref{SO3dexpInv2}) or (\ref%
{SO3dexpInv3}) calculated using (\ref{DD}) in Lemma \ref{LemmaSE3dexp}.
\end{proof}

It was shown in \cite[p. 77]{SeligBook} that the matrix of the
right-trivialized differential and its inverse can be expressed directly in
terms of the adjoint operator matrix on $se\left( 3\right) $ in (\ref{adSE3}%
), which satisfies $\mathbf{ad}_{\mathbf{X}}^{6}=-\left\Vert \mathbf{x}%
\right\Vert ^{4}\mathbf{ad}_{\mathbf{X}}^{2}-2\left\Vert \mathbf{x}%
\right\Vert ^{2}\mathbf{ad}_{\mathbf{X}}^{4}$. In terms of the parameters (%
\ref{abc}), they are 
\begin{align}
\mathbf{dexp}_{\mathbf{X}}& =\mathbf{I}+\tfrac{1}{2\left\Vert \mathbf{x}%
\right\Vert ^{2}}\left( 4-\left\Vert \mathbf{x}\right\Vert \sin \left\Vert 
\mathbf{x}\right\Vert -4\cos \left\Vert \mathbf{x}\right\Vert \right) 
\mathbf{ad}_{\mathbf{X}}+\tfrac{1}{2\left\Vert \mathbf{x}\right\Vert ^{3}}%
\left( 4\left\Vert \mathbf{x}\right\Vert -5\sin \left\Vert \mathbf{x}%
\right\Vert +\left\Vert \mathbf{x}\right\Vert \cos \left\Vert \mathbf{x}%
\right\Vert \right) \mathbf{ad}_{\mathbf{X}}^{2}  \notag \\
& +\tfrac{1}{2\left\Vert \mathbf{x}\right\Vert ^{4}}\left( 2-\left\Vert 
\mathbf{x}\right\Vert \sin \left\Vert \mathbf{x}\right\Vert -2\cos
\left\Vert \mathbf{x}\right\Vert \right) \mathbf{ad}_{\mathbf{X}}^{3}+\tfrac{%
1}{2\left\Vert \mathbf{x}\right\Vert ^{5}}\left( 2\left\Vert \mathbf{x}%
\right\Vert -3\sin \left\Vert \mathbf{x}\right\Vert +\left\Vert \mathbf{x}%
\right\Vert \cos \left\Vert \mathbf{x}\right\Vert \right) \mathbf{ad}_{%
\mathbf{X}}^{4}  \label{dexpSE3ad} \\
& =\mathbf{I}+\left( \beta -\tfrac{\alpha }{2}\right) \mathbf{ad}_{\mathbf{X}%
}+\frac{1}{2}\left( 5\delta -\tfrac{\beta }{2}\right) \mathbf{ad}_{\mathbf{X}%
}^{2}+\tfrac{1}{2\left\Vert \mathbf{x}\right\Vert ^{2}}\left( \beta -\alpha
\right) \mathbf{ad}_{\mathbf{X}}^{3}+\tfrac{1}{2\left\Vert \mathbf{x}%
\right\Vert ^{2}}\left( 3\delta -\tfrac{\beta }{2}\right) \mathbf{ad}_{%
\mathbf{X}}^{4}  \label{dexpSE3ad2} \\
\mathbf{dexp}_{\mathbf{X}}^{-1}& =\mathbf{I}-\frac{1}{2}\mathbf{ad}_{\mathbf{%
X}}+\left( \tfrac{2}{\left\Vert \mathbf{x}\right\Vert ^{2}}+\tfrac{%
\left\Vert \mathbf{x}\right\Vert +3\sin \left\Vert \mathbf{x}\right\Vert }{%
4\left\Vert \mathbf{x}\right\Vert \left( \cos \left\Vert \mathbf{x}%
\right\Vert -1\right) }\right) \mathbf{ad}_{\mathbf{X}}^{2}+\left( \tfrac{1}{%
\left\Vert \mathbf{x}\right\Vert ^{4}}+\tfrac{\left\Vert \mathbf{x}%
\right\Vert +\sin \left\Vert \mathbf{x}\right\Vert }{4\left\Vert \mathbf{x}%
\right\Vert ^{3}\left( \cos \left\Vert \mathbf{x}\right\Vert -1\right) }%
\right) \mathbf{ad}_{\mathbf{X}}^{4}  \label{dexpInvSE3ad} \\
& =\mathbf{I}-\frac{1}{2}\mathbf{ad}_{\mathbf{X}}+\tfrac{1}{\left\Vert 
\mathbf{x}\right\Vert ^{2}}\left( 2-\tfrac{1+3\alpha }{2\beta }\right) 
\mathbf{ad}_{\mathbf{X}}^{2}+\tfrac{1}{\left\Vert \mathbf{x}\right\Vert ^{4}}%
\left( 1-\tfrac{1+\alpha }{2\beta }\right) \mathbf{ad}_{\mathbf{X}}^{4}.
\label{dexpInvSE3ad2}
\end{align}%
The expression (\ref{dexpInvSE3ad}) was already presented in \cite%
{BulloMurray1995} with negative sign so to express the left-trivialized
differential. Using different abbreviations, these expressions were reported
in \cite{BottassoBorri1998}, where it arose naturally using the adjoint
representation (see section \ref{secAdjSE3}). In contrast to (\ref%
{diffDexpInvSO3}), the expression (\ref{dexpInvSE3ad2}) can be evaluated
numerically stable without separating the case when $\mathbf{x}=\mathbf{0}$.
If efficiently implemented, it provides a computationally efficient
alternative to (\ref{DexpInvSE3}).

\subsection{Directional derivative of the right-trivialized differential}

The computational procedure of the (semi)implicit generalized-$\alpha $ Lie
group method \cite%
{ArnoldBrulsCardona2015,ArnoldHante2017,BrulsCardona2010,BrulsCardonaArnold2012}
employs an iteration matrix (called the tangent matrix), which involves the
directional derivative of the trivialized differential of the respective
coordinate map. This is necessary also when using other implicit integration
schemes.

To simplify notation, denote the directional derivative (\ref{diffDexpSO33})
of the matrix $\mathbf{dexp}_{\mathbf{x}}$ in (\ref{dexpSO35}) with $\mathbf{%
Ddexp}\left( \mathbf{X}\right) :=(\mathrm{D}_{\mathbf{x}}\mathbf{dexp})%
%TCIMACRO{\TeXButton{-0.5ex}{\hspace{-0.5ex}}}%
%BeginExpansion
\hspace{-0.5ex}%
%EndExpansion
\left( \mathbf{y}\right) $, and the derivative (\ref{diffDexpInvSO3}) of $%
\mathbf{dexp}_{\mathbf{x}}^{-1}$ in (\ref{SO3dexpInv3}) with $\mathbf{Ddexp}%
^{-1}\left( \mathbf{X}\right) :=(\mathrm{D}_{\mathbf{x}}\mathbf{dexp}^{-1})%
%TCIMACRO{\TeXButton{-0.5ex}{\hspace{-0.5ex}}}%
%BeginExpansion
\hspace{-0.5ex}%
%EndExpansion
\left( \mathbf{y}\right) $, where $\mathbf{X}=\left( \mathbf{x},\mathbf{y}%
\right) $. The second directional derivative of $\mathbf{dexp}$ is then
written as $(\mathrm{D}_{\mathbf{X}}\mathbf{Ddexp})%
%TCIMACRO{\TeXButton{-0.5ex}{\hspace{-0.5ex}}}%
%BeginExpansion
\hspace{-0.5ex}%
%EndExpansion
\left( \mathbf{U}\right) $ with $\mathbf{U}=\left( \mathbf{u},\mathbf{v}%
\right) $, and similarly for its inverse.

\begin{lemma}
The directional derivative of the matrix $\mathbf{dexp}$ in (\ref{dexpSE31})
is%
\begin{equation}
(\mathrm{D}_{\mathbf{X}}\mathbf{dexp})%
%TCIMACRO{\TeXButton{-0.5ex}{\hspace{-0.5ex}}}%
%BeginExpansion
\hspace{-0.5ex}%
%EndExpansion
\left( \mathbf{U}\right) =\left( 
\begin{array}{cc}
\mathrm{D}_{\mathbf{x}}\mathbf{dexp})%
%TCIMACRO{\TeXButton{-0.5ex}{\hspace{-0.5ex}}}%
%BeginExpansion
\hspace{-0.5ex}%
%EndExpansion
\left( \mathbf{u}\right) & \mathbf{0} \\ 
(\mathrm{D}_{\mathbf{X}}\mathbf{Ddexp})%
%TCIMACRO{\TeXButton{-0.5ex}{\hspace{-0.5ex}}}%
%BeginExpansion
\hspace{-0.5ex}%
%EndExpansion
\left( \mathbf{U}\right) & (\mathrm{D}_{\mathbf{x}}\mathbf{dexp})%
%TCIMACRO{\TeXButton{-0.5ex}{\hspace{-0.5ex}}}%
%BeginExpansion
\hspace{-0.5ex}%
%EndExpansion
\left( \mathbf{u}\right)%
\end{array}%
\right)  \label{diffDexpSE3}
\end{equation}%
where $\mathbf{U}=\left( \mathbf{u},\mathbf{v}\right) $, and the directional
derivative of matrix $\mathbf{Ddexp}$ possesses the explicitly form%
\begin{align}
(\mathrm{D}_{\mathbf{X}}\mathbf{Ddexp})%
%TCIMACRO{\TeXButton{-0.5ex}{\hspace{-0.5ex}}}%
%BeginExpansion
\hspace{-0.5ex}%
%EndExpansion
\left( \mathbf{U}\right) =~& \tfrac{\beta }{2}\tilde{\mathbf{v}}+\tfrac{%
\alpha -\beta }{\left\Vert \mathbf{x}\right\Vert ^{2}}\left( (\mathbf{x}^{T}%
\mathbf{u})\tilde{\mathbf{y}}+(\mathbf{x}^{T}\mathbf{v+y}^{T}\mathbf{u})%
\tilde{\mathbf{x}}+(\mathbf{x}^{T}\mathbf{y})\tilde{\mathbf{u}}\right)
\label{diff2dexpSO31} \\
& -\tfrac{\beta }{2\left\Vert \mathbf{x}\right\Vert ^{2}}(\mathbf{x}^{T}%
\mathbf{u})(\mathbf{x}^{T}\mathbf{y})\tilde{\mathbf{x}}+\delta \left( \tilde{%
\mathbf{x}}\tilde{\mathbf{v}}+\tilde{\mathbf{v}}\tilde{\mathbf{x}}+\tilde{%
\mathbf{y}}\tilde{\mathbf{u}}+\tilde{\mathbf{u}}\tilde{\mathbf{y}}\right) 
\notag \\
& +\tfrac{\beta /2-3\delta }{\left\Vert \mathbf{x}\right\Vert ^{2}}\left( (%
\mathbf{x}^{T}\mathbf{y})(\tilde{\mathbf{x}}\tilde{\mathbf{u}}+\tilde{%
\mathbf{u}}\tilde{\mathbf{x}})+(\mathbf{x}^{T}\mathbf{u})(\tilde{\mathbf{x}}%
\tilde{\mathbf{y}}+\tilde{\mathbf{y}}\tilde{\mathbf{x}})+(\mathbf{x}^{T}%
\mathbf{v}+\mathbf{y}^{T}\mathbf{u})\tilde{\mathbf{x}}^{2}\right)  \notag \\
& +\tfrac{1}{\left\Vert \mathbf{x}\right\Vert ^{4}}(\mathbf{x}^{T}\mathbf{y}%
)(\mathbf{x}^{T}\mathbf{u})\left( (1-5\alpha +4\beta )\tilde{\mathbf{x}}%
+(\alpha -\tfrac{7}{2}\beta +15\delta )\tilde{\mathbf{x}}^{2}\right)  \notag
\\
=~& \tfrac{\beta }{2}\tilde{\mathbf{v}}+\tfrac{\alpha -\beta }{\left\Vert 
\mathbf{x}\right\Vert ^{2}}(\mathbf{x}^{T}\mathbf{u})\tilde{\mathbf{y}} 
\label{diff2dexpSO32} \\
& 
%TCIMACRO{\TeXButton{red}{\color{red}}}%
%BeginExpansion
\color{red}%
%EndExpansion
+\tfrac{\beta /2-3\delta }{\left\Vert \mathbf{x}\right\Vert ^{2}}\left( (%
\mathbf{x}^{T}\mathbf{u})(\tilde{\mathbf{x}}\tilde{\mathbf{y}}+\tilde{%
\mathbf{y}}\tilde{\mathbf{x}})+(\mathbf{x}^{T}\mathbf{y})(\tilde{\mathbf{x}}%
\tilde{\mathbf{u}}+\tilde{\mathbf{u}}\tilde{\mathbf{x}})\right) +\delta
\left( \tilde{\mathbf{x}}\tilde{\mathbf{v}}+\tilde{\mathbf{v}}\tilde{\mathbf{%
x}}+\tilde{\mathbf{y}}\tilde{\mathbf{u}}+\tilde{\mathbf{u}}\tilde{\mathbf{y}}%
\right)   \notag \\
& +\tfrac{1}{\left\Vert \mathbf{x}\right\Vert ^{2}}\left( \left( \alpha
-\beta \right) (\mathbf{x}^{T}\mathbf{v+y}^{T}\mathbf{u})-\tfrac{1}{2}\beta (%
\mathbf{x}^{T}\mathbf{u})(\mathbf{x}^{T}\mathbf{y})+\tfrac{1-5\alpha +4\beta 
}{\left\Vert \mathbf{x}\right\Vert ^{2}}(\mathbf{x}^{T}\mathbf{y})(\mathbf{x}%
^{T}\mathbf{u})\right) \tilde{\mathbf{x}}  \notag \\
& +\tfrac{1}{\left\Vert \mathbf{x}\right\Vert ^{2}}\left( (\tfrac{1}{2}\beta
-3\delta )(\mathbf{x}^{T}\mathbf{v+y}^{T}\mathbf{u})+\tfrac{1}{\left\Vert 
\mathbf{x}\right\Vert ^{2}}(\alpha -\tfrac{7}{2}\beta +15\delta )(\mathbf{x}%
^{T}\mathbf{y})(\mathbf{x}^{T}\mathbf{u})\right) \tilde{\mathbf{x}}^{2} 
\notag \\
& 
%TCIMACRO{\TeXButton{red}{\color{red}}}%
%BeginExpansion
\color{red}%
%EndExpansion
+\tfrac{\alpha -\beta }{\left\Vert \mathbf{x}\right\Vert ^{2}}( \mathbf{%
x}^{T}\mathbf{y}) \tilde{\mathbf{u}}.  \notag
\end{align}
\end{lemma}

\begin{proof}
The derivatives of the parameters (\ref{abc}) are readily found with (\ref%
{DD}) to be%
\begin{equation}
(\mathrm{D}_{\mathbf{x}}\alpha )%
%TCIMACRO{\TeXButton{-0.5ex}{\hspace{-0.5ex}}}%
%BeginExpansion
\hspace{-0.5ex}%
%EndExpansion
\left( \mathbf{u}\right) =\left( \delta -\tfrac{1}{2}\beta \right) \mathbf{x}%
^{T}\mathbf{u},(\mathrm{D}_{\mathbf{x}}\beta )%
%TCIMACRO{\TeXButton{-0.5ex}{\hspace{-0.5ex}}}%
%BeginExpansion
\hspace{-0.5ex}%
%EndExpansion
\left( \mathbf{u}\right) =\tfrac{2\mathbf{x}^{T}\mathbf{u}}{\left\Vert 
\mathbf{x}\right\Vert ^{2}}\left( \alpha -\beta \right) ,(\mathrm{D}_{%
\mathbf{x}}\delta )%
%TCIMACRO{\TeXButton{-0.5ex}{\hspace{-0.5ex}}}%
%BeginExpansion
\hspace{-0.5ex}%
%EndExpansion
\left( \mathbf{u}\right) =\tfrac{\mathbf{x}^{T}\mathbf{u}}{\left\Vert 
\mathbf{x}\right\Vert ^{2}}\left( \tfrac{1}{2}\beta -3\delta \right)
\label{DD2}
\end{equation}%
A straightforward manipulation using (\ref{DD2}) yields the directional
derivative of (\ref{diffDexpSO33}) in the form (\ref{diff2dexpSO31}), and
rearranging yields (\ref{diff2dexpSO32}). The statement follows with the
matrix of $\mathbf{dexp}_{\mathbf{X}}$ in (\ref{dexpSE31}).
\end{proof}

The directional derivative (\ref{diffDexpSE3}) can be evaluated along with
one of the expressions in (\ref{diffDexpSO34})-(\ref{diffDexpSO33}). A
slightly different expression for the second directional derivative of the
left-trivialized differential was presented in \cite{SonnevillePhD}.

\begin{lemma}
The directional derivative of the matrix $\mathbf{dexp}^{-1}$ in (\ref%
{DexpInvSE3}) is%
\begin{equation}
(\mathrm{D}_{\mathbf{X}}\mathbf{dexp}^{-1})%
%TCIMACRO{\TeXButton{-0.5ex}{\hspace{-0.5ex}}}%
%BeginExpansion
\hspace{-0.5ex}%
%EndExpansion
\left( \mathbf{U}\right) =\left( 
\begin{array}{cc}
(\mathrm{D}_{\mathbf{x}}\mathbf{dexp}^{-1})%
%TCIMACRO{\TeXButton{-0.5ex}{\hspace{-0.5ex}}}%
%BeginExpansion
\hspace{-0.5ex}%
%EndExpansion
\left( \mathbf{u}\right)  & \mathbf{0} \\ 
(\mathrm{D}_{\mathbf{X}}\mathbf{Ddexp}^{-1})%
%TCIMACRO{\TeXButton{-0.5ex}{\hspace{-0.5ex}}}%
%BeginExpansion
\hspace{-0.5ex}%
%EndExpansion
\left( \mathbf{U}\right)  & (\mathrm{D}_{\mathbf{x}}\mathbf{dexp}^{-1})%
%TCIMACRO{\TeXButton{-0.5ex}{\hspace{-0.5ex}}}%
%BeginExpansion
\hspace{-0.5ex}%
%EndExpansion
\left( \mathbf{u}\right) 
\end{array}%
\right)   \label{diffDexpInvSE3}
\end{equation}%
where $\mathbf{U}=\left( \mathbf{u},\mathbf{v}\right) $, and the directional
derivative of matrix $\mathbf{Ddexp}^{-1}$ possesses the explicit form%
\begin{align}
(\mathrm{D}_{\mathbf{X}}\mathbf{Ddexp}^{-1})%
%TCIMACRO{\TeXButton{-0.5ex}{\hspace{-0.5ex}}}%
%BeginExpansion
\hspace{-0.5ex}%
%EndExpansion
\left( \mathbf{U}\right) & =-\tfrac{1}{2}\tilde{\mathbf{v}}+%
%TCIMACRO{\TeXButton{red}{\color{red}}}%
%BeginExpansion
\color{red}%
%EndExpansion
\tfrac{1}{\left\Vert \mathbf{x}\right\Vert ^{2}}\left( (1-\gamma )\left( 
\tilde{\mathbf{x}}\tilde{\mathbf{v}}+\tilde{\mathbf{v}}\tilde{\mathbf{x}}+%
\tilde{\mathbf{y}}\tilde{\mathbf{u}}+\tilde{\mathbf{u}}\tilde{\mathbf{y}}%
\right) +\tfrac{1}{4}(\mathbf{x}^{T}\mathbf{u})(\tilde{\mathbf{x}}\tilde{%
\mathbf{y}}+\tilde{\mathbf{y}}\tilde{\mathbf{x}})\right)
\label{diff2dexpInvSO31} \\
& 
%TCIMACRO{\TeXButton{red}{\color{red}}}%
%BeginExpansion
\color{red}%
%EndExpansion
+\tfrac{1}{\left\Vert \mathbf{x}\right\Vert ^{4}}%
%TCIMACRO{\TeXButton{Big}{\Big}}%
%BeginExpansion
\Big%
%EndExpansion
(\left( \gamma -1\right) (2+\gamma )(\mathbf{x}^{T}\mathbf{u})(\tilde{%
\mathbf{x}}\tilde{\mathbf{y}}+\tilde{\mathbf{y}}\tilde{\mathbf{x}})  \notag
\\
& 
%TCIMACRO{\TeXButton{red}{\color{red}}}%
%BeginExpansion
\color{red}%
%EndExpansion
+(\gamma +\tfrac{1}{\beta }-2)\left( (\mathbf{x}^{T}\mathbf{y})(\tilde{%
\mathbf{x}}\tilde{\mathbf{u}}+\tilde{\mathbf{u}}\tilde{\mathbf{x}})+(\mathbf{%
x}^{T}\mathbf{v}+\mathbf{y}^{T}\mathbf{u})\tilde{\mathbf{x}}^{2}\right) -%
\tfrac{1}{4}(\mathbf{x}^{T}\mathbf{y})(\mathbf{x}^{T}\mathbf{u})\tilde{%
\mathbf{x}}^{2}%
%TCIMACRO{\TeXButton{Big}{\Big}}%
%BeginExpansion
\Big%
%EndExpansion
)  \notag \\
& +\tfrac{1}{\left\Vert \mathbf{x}\right\Vert ^{6}}(\mathbf{x}^{T}\mathbf{y}%
)(\mathbf{x}^{T}\mathbf{u})\left( 8-3\gamma -\gamma ^{2}-\tfrac{2}{\beta ^{2}%
}\left( \alpha +\beta \right) \right) \tilde{\mathbf{x}}^{2}  \notag \\
=& 
%TCIMACRO{\TeXButton{red}{\red}}%
%BeginExpansion
\red%
%EndExpansion
-\tfrac{1}{2}\tilde{\mathbf{v}}+\tfrac{1}{\left\Vert \mathbf{x}\right\Vert
^{2}}\left( \tfrac{1}{4}(\mathbf{x}^{T}\mathbf{u})(\tilde{\mathbf{x}}\tilde{%
\mathbf{y}}+\tilde{\mathbf{y}}\tilde{\mathbf{x}})+\left( 1-\gamma \right)
\left( \tilde{\mathbf{x}}\tilde{\mathbf{v}}+\tilde{\mathbf{v}}\tilde{\mathbf{%
x}}+\tilde{\mathbf{y}}\tilde{\mathbf{u}}+\tilde{\mathbf{u}}\tilde{\mathbf{y}}%
\right) \right)  \label{diff2dexpInvSO32} \\
& 
%TCIMACRO{\TeXButton{red}{\red}}%
%BeginExpansion
\red%
%EndExpansion
+\tfrac{1}{\left\Vert \mathbf{x}\right\Vert ^{4}}%
%TCIMACRO{\TeXButton{Big}{\Big}}%
%BeginExpansion
\Big%
%EndExpansion
((\mathbf{x}^{T}\mathbf{u})\left( \gamma \left( 1+\gamma \right) -2\right) (%
\tilde{\mathbf{x}}\tilde{\mathbf{y}}+\tilde{\mathbf{y}}\tilde{\mathbf{x}})+(%
\tfrac{1}{\beta }+\gamma -2)(\mathbf{x}^{T}\mathbf{y})(\tilde{\mathbf{x}}%
\tilde{\mathbf{u}}+\tilde{\mathbf{u}}\tilde{\mathbf{x}})  \notag \\
& 
%TCIMACRO{\TeXButton{red}{\red}}%
%BeginExpansion
\red%
%EndExpansion
%TCIMACRO{\TeXButton{TeX field}{\hspace{8ex}}}%
%BeginExpansion
\hspace{8ex}%
%EndExpansion
+%
%TCIMACRO{\TeXButton{big}{\big}}%
%BeginExpansion
\big%
%EndExpansion
((\tfrac{1}{\beta }+\gamma -2)(\mathbf{x}^{T}\mathbf{v}+\mathbf{y}^{T}%
\mathbf{u})-\tfrac{1}{4}(\mathbf{x}^{T}\mathbf{y})(\mathbf{x}^{T}\mathbf{u})%
%TCIMACRO{\TeXButton{big}{\big}}%
%BeginExpansion
\big%
%EndExpansion
)\tilde{\mathbf{x}}^{2}%
%TCIMACRO{\TeXButton{Big}{\Big}}%
%BeginExpansion
\Big%
%EndExpansion
)  \notag \\
& 
%TCIMACRO{\TeXButton{red}{\red}}%
%BeginExpansion
\red%
%EndExpansion
+\tfrac{1}{\left\Vert \mathbf{x}\right\Vert ^{6}}(\mathbf{x}^{T}\mathbf{y})(%
\mathbf{x}^{T}\mathbf{u})\left( 8-\gamma \left( 3+\gamma \right) -\tfrac{2}{%
\beta }\left( 1+\gamma \right) \right) \tilde{\mathbf{x}}^{2}.  \notag
\end{align}
\end{lemma}

\begin{proof}
The directional derivative of (2.40) involves the derivative of parameter $%
\gamma $ 
\begin{equation}
(\mathrm{D}_{\mathbf{x}}\gamma )%
%TCIMACRO{\TeXButton{-0.5ex}{\hspace{-0.5ex}}}%
%BeginExpansion
\hspace{-0.5ex}%
%EndExpansion
\left( \mathbf{u}\right) =\tfrac{\mathbf{x}^{T}\mathbf{u}}{\left\Vert 
\mathbf{x}\right\Vert ^{2}}\left( \gamma -\gamma ^{2}-\tfrac{\left\Vert 
\mathbf{x}\right\Vert ^{2}}{4}\right) .
\end{equation}%
The derivative of the coefficient of $\tilde{\mathbf{x}}^{2}$ in (2.40) is
then found to be, along with (2.50),%
\begin{align}
\left( \mathrm{D}_{\mathbf{X}}(\tfrac{\mathbf{x}^{T}\mathbf{y}}{\left\Vert 
\mathbf{x}\right\Vert ^{4}}(\tfrac{1}{\beta }+\gamma -2))\right) 
%TCIMACRO{\TeXButton{-0.5ex}{\hspace{-0.5ex}}}%
%BeginExpansion
\hspace{-0.5ex}%
%EndExpansion
\left( \mathbf{U}\right) =& \ 
%TCIMACRO{\TeXButton{red}{\color{red}} }%
%BeginExpansion
\color{red}
%EndExpansion
\tfrac{1}{\left\Vert \mathbf{x}\right\Vert ^{4}}(\tfrac{1}{\beta }+\gamma
-2)(\mathbf{x}^{T}\mathbf{v}+\mathbf{y}^{T}\mathbf{u})-\tfrac{1}{4}(\mathbf{x%
}^{T}\mathbf{y})(\mathbf{x}^{T}\mathbf{u}))  \notag \\
& 
%TCIMACRO{\TeXButton{red}{\color{red}} }%
%BeginExpansion
\color{red}
%EndExpansion
+\tfrac{(\mathbf{x}^{T}\mathbf{y})(\mathbf{x}^{T}\mathbf{u})}{\left\Vert 
\mathbf{x}\right\Vert ^{6}}\left( 8-\gamma \left( 3+\gamma \right) -\tfrac{2%
}{\beta }\left( 1+\gamma \right) \right) .
\end{align}%
The expressions and (\ref{diff2dexpInvSO31}) and (\ref{diff2dexpInvSO32})
are obtained after some algebraic manipulation.
\end{proof}

An equivalent expression for the directional derivative of the matrix of the
left-trivialized differential, i.e. with negative argument $\mathbf{X}$) was
derived in \cite[appendix A.1]{SonnevillePhD}. The directional derivative is
thus available in closed form as (\ref{diffDexpInvSE3}) along with (\ref%
{diffDexpInvSO3}).

The above relations for the derivative of the matrices representing the
directional derivatives of the right-trivialized differential and its
inverse are crucial for numerical integration using implicit Lie group
integration methods, such as the generalized-$\alpha $ scheme. They can be
implemented so to cope with $\left\Vert \mathbf{x}\right\Vert =0$. An
alternative formulation could be obtained using the expressions (\ref%
{dexpSE3ad2}) and (\ref{dexpInvSE3ad2}), respectively.

\subsection{Adjoint Representation -- The 'Configuration Tensor'%
%TCIMACRO{\TeXButton{secAdjSE3}{\label{secAdjSE3}}}%
%BeginExpansion
\label{secAdjSE3}%
%EndExpansion
}

Representing frame transformations by matrices $\mathbf{C}\in SE\left(
3\right) $, which describe transformations of homogenous point coordinates,
is useful to compute relative configurations of rigid bodies by means of
matrix multiplication. This is not relevant when $\mathbf{C}$ describes the
absolute configuration of a rigid body or the displacement field of Cosserat
beam, for instances. Moreover, storing rotation matrix and position vector
in $\mathbf{C}$ is merely a means of bookkeeping. In this context a more
relevant operation is the frame transformation of twists, which is described
by the adjoint operator $\mathrm{Ad}_{\mathbf{C}}:se\left( 3\right)
\rightarrow se\left( 3\right) $. In vector representation of twists, this is
expressed by the matrix \cite{SeligBook,MurrayBook}%
\begin{equation}
\mathbf{Ad}_{\mathbf{C}}=\left( 
\begin{array}{cc}
\mathbf{R} & \mathbf{0} \\ 
\tilde{\mathbf{r}}\mathbf{R} & \mathbf{R}%
\end{array}%
\right) \in \mathrm{Ad}\left( SE\left( 3\right) \right) .  \label{AdSE3}
\end{equation}%
This \emph{adjoint representation} is used in the geometrically exact beam
formulation by Borri et al. \cite%
{BorriMelloAtluri1990,Borri2001a,Borri2001b,BorriBottassoTrainelli2003} and
Bauchau \cite{BauchauBook2010}. In this context, the adjoint transformation
matrix $\mathbf{Ad}_{\mathbf{C}}$ is referred to as the \emph{configuration
tensor} \cite{BauchauChoi2003,BauchauBook2010}, and used to describe frame
transformations, where $\mathrm{Ad}\left( SE\left( 3\right) \right) $ was
denoted $SR\left( 6\right) $. The expression (\ref{Poisson}) for the twists,
and (\ref{deform}) for the deformation tensor defining the strain field, are
then replaced by%
\begin{align}
\dot{\mathbf{Ad}}_{\mathbf{C}}& =\mathbf{ad}_{{\mathbf{V}}{^{\mathrm{s}}}}%
\mathbf{Ad}_{\mathbf{C}},\ \dot{\mathbf{Ad}_{\mathbf{C}}}=\mathbf{Ad}_{%
\mathbf{C}}\mathbf{ad}_{{\mathbf{V}}{^{\mathrm{b}}}} \\
\mathbf{Ad}_{\mathbf{C}}^{\prime }& =\mathbf{ad}_{{\bm{\chi}}{^{\mathrm{s}}}}%
\mathbf{Ad}_{\mathbf{C}},\ \mathbf{Ad}_{\mathbf{C}}^{\prime }=\mathbf{Ad}_{%
\mathbf{C}}\mathbf{ad}_{{\bm{\chi}}{^{\mathrm{b}}}}
\end{align}%
where $\mathbf{ad}_{{\mathbf{V}}}$ and $\mathbf{ad}_{{\bm{\chi}}}$ serve to
represent rigid body twists and strain measures, with $\mathbf{ad}$ in (\ref%
{adSE3}).

The adjoint representation is canonically parameterized in terms of
instantaneous screw coordinates $\mathbf{X}\in {\mathbb{R}}^{6}$ using the
relation%
\begin{equation}
\mathbf{Ad}_{\exp \hat{\mathbf{X}}}=\exp \mathbf{ad}_{\mathbf{X}}
\label{Adexp}
\end{equation}%
where $\exp $ is the exponential map of the adjoint representation. Clearly
the canonical (screw) coordinates $\mathbf{X}\in {\mathbb{R}}^{6}$ in the
exponential map (\ref{expSE3}) and (\ref{Adexp}) are identical. An important
consequence is that the local reconstruction equations (\ref{RecExp}) remain
valid, and hence the closed form expressions (\ref{dexpSE31}) and (\ref%
{DexpInvSE3}) apply. Closed form expressions for the exponential (\ref{Adexp}%
) and its derivatives were reported in \cite{BorriTrainelliBottasso2000}.
This will not be considered here as recent modeling approaches are using the 
$SE\left( 3\right) $ representation of Euclidean motions \cite%
{SonnevilleCardonaBruls2014,SonnevillePhD}. Moreover, the adjoint
representation (\ref{AdSE3}) is redundant and does not directly yield the
displacement vector $\mathbf{r}$.

For completeness, the following relations, which were central for developing
integration schemes using the base-pole formulation \cite%
{Borri2001a,Borri2001b}, and reported without proof, are summarized.

\begin{lemma}
Let $\mathbf{C}=\exp \hat{\mathbf{X}}\in SE\left( 3\right) $, the adjoint
operator matrix (\ref{AdSE3}) and right-trivialized differential of $\exp
:se\left( 3\right) \rightarrow SE\left( 3\right) $ satisfy the relations%
\begin{align}
\mathbf{Ad}_{\mathbf{C}} &=\mathbf{dexp}_{-\hat{\mathbf{X}}}^{-1}\mathbf{dexp%
}_{\hat{\mathbf{X}}}  \label{SE3dexpdexp1} \\
&=\mathbf{dexp}_{\hat{\mathbf{X}}}\mathbf{dexp}_{-\hat{\mathbf{X}}}^{-1}
\label{SE3dexpdexp2} \\
&=\mathbf{I}+\mathbf{ad}_{\hat{\mathbf{X}}}\mathbf{dexp}_{\hat{\mathbf{X}}}
\label{SE3dexpdexp3} \\
&=\mathbf{I}+\mathbf{dexp}_{\hat{\mathbf{X}}}\mathbf{ad}_{\hat{\mathbf{X}}}.
\label{SE3dexpdexp4}
\end{align}
\end{lemma}

\begin{proof}
The proof is similar to that of lemma \ref{LemmaSO3dexp}. From (\ref{diff1})
and (\ref{diff2}) follows $\mathrm{dexp}_{\hat{\mathbf{X}}}(\hat{\mathbf{Y}}%
)=\mathbf{C}\mathrm{dexp}_{-\hat{\mathbf{X}}}(\hat{\mathbf{Y}})\mathbf{C}%
^{-1}=\mathrm{Ad}_{\mathbf{C}}(\mathrm{dexp}_{-\hat{\mathbf{X}}}(\hat{%
\mathbf{Y}}))$, and hence the vector representation of $se\left( 3\right) $
yields $\mathbf{dexp}_{\mathbf{X}}\mathbf{Y}=\mathbf{Ad}_{\mathbf{C}}\mathbf{%
dexp}_{-\mathbf{X}}\mathbf{Y}$ and thus (\ref{SE3dexpdexp1}). The explicit
form (\ref{dexpInvSE3ad}) shows that $\mathbf{ad}_{\mathbf{X}}=\mathbf{dexp}%
_{-\mathbf{X}}^{-1}-\mathbf{dexp}_{\mathbf{X}}^{-1}$, and hence $\mathbf{I}+%
\mathbf{dexp}_{\mathbf{X}}\mathbf{ad}_{\mathbf{X}}=\mathbf{dexp}_{\mathbf{X}}%
\mathbf{dexp}_{-\mathbf{X}}^{-1}$ in (\ref{SE3dexpdexp4}). From (\ref%
{SO3dexpInv1}) follows that $\mathbf{ad}_{\mathbf{X}}\mathbf{dexp}_{\mathbf{X%
}}=\mathbf{dexp}_{\mathbf{X}}\mathbf{ad}_{\mathbf{X}}$, which proves (\ref%
{SE3dexpdexp3}).
\end{proof}

The relation (\ref{SE3dexpdexp1})-(\ref{SE3dexpdexp4}) were exploited in 
\cite{Borri2001a,Borri2001b} to construct invariant conserving, respectively
dissipating, integration schemes for multibody systems with flexible
members. Notice that (\ref{SE3dexpdexp1})-(\ref{SE3dexpdexp4}) hold true for
general matrix Lie groups. The relations (\ref{dexpdexpSO31})-(\ref%
{dexpdexpSO34}) are special cases with $\mathbf{ad}_{\tilde{\mathbf{x}}}=%
\tilde{\mathbf{x}},\mathbf{Ad}_{\mathbf{R}}=\mathbf{R}$.

\section{Motion Parameterization via the Cayley Map%
%TCIMACRO{\TeXButton{secCay}{\label{secCay}}}%
%BeginExpansion
\label{secCay}%
%EndExpansion
}

The vectorial parameterizations of motion admit an algebraic description
without transcendental functions. Parameterizations in terms of
non-redundant parameters provide a computationally efficient alternative to
canonical coordinates within Lie group integration schemes. The simplest of
such are the Rodrigues parameters. The corresponding coordinate maps are
obtained via the Cayley map.

The Cayley map on a quadratic Lie group provides an approximation of the
exponential map \cite{HairerLubichWanner2006}. Aiming at computationally
efficient Lie group integration schemes, the Cayley map has been widely used
for general systems \cite%
{DieleLopezPeluso1998,DieleLopezPoliti1998,LewisSimo1994,BlanesCasasRos2002}%
, and for the rigid body motion described on $SO\left( 3\right) $ in
particular, e.g. \cite{CelledoniOwren2003,IserlesMuntheKaasNrsettZanna2000}.
A symplectic integration scheme for rotating rigid bodies, where the Cayley
map is used, was presented in \cite{LewisSimo1994}. The Cayley map for the
adjoint representation of $SE\left( 3\right) $ along with the base-pol
formulation was used to derive integration schemes for geometrically exact
beams \cite{BorriBottassoTrainelli2003} that preserve certain invariants.
The necessary explicit relations for the Cayley map and its differentials
are scattered in the literature partly without proof. Moreover, to the
author's knowledge, the differential of the Cayley map on $SE\left( 3\right) 
$ and its directional derivative are not present in the literature. They
will be presented in this section.

The Cayley map $\mathrm{cay}:\mathfrak{g}\rightarrow G$ on a quadratic
matrix Lie group can be defined as%
\begin{align}
\mathrm{cay}\left( \tilde{\mathbf{x}}\right) & =\left( \mathbf{I}-\tilde{%
\mathbf{x}}\right) ^{-1}\left( \mathbf{I}+\tilde{\mathbf{x}}\right)
\label{Cay1} \\
& =\left( \mathbf{I}+\tilde{\mathbf{x}}\right) \left( \mathbf{I}-\tilde{%
\mathbf{x}}\right) ^{-1}  \label{Cay4} \\
& =\mathbf{I}+2(\tilde{\mathbf{x}}+\tilde{\mathbf{x}}^{2}+\tilde{\mathbf{x}}%
^{3}+\ldots ).  \label{Cay3}
\end{align}%
The components of $\mathbf{x}\in {\mathbb{R}}^{n}$ serve as local
(non-canonical) coordinates around the identity in $G$. The directional
derivative is 
\begin{align}
\left( \mathrm{D}_{\tilde{\mathbf{x}}}\mathrm{cay}\right) \left( \tilde{%
\mathbf{y}}\right) & :=\frac{d}{dt}\mathrm{cay}(\hat{\mathbf{X}}+t\hat{%
\mathbf{Y}})|_{t=0}  \notag \\
& =\left( \mathbf{I}-\tilde{\mathbf{x}}\right) ^{-1}\tilde{\mathbf{y}}\left( 
\mathbf{I}+\mathrm{cay}\left( \tilde{\mathbf{x}}\right) \right) =2\left( 
\mathbf{I}-\tilde{\mathbf{x}}\right) ^{-1}\tilde{\mathbf{y}}\left( \mathbf{I}%
-\tilde{\mathbf{x}}\right) ^{-1}  \label{Cay2}
\end{align}%
and thus, with (\ref{Cay1}), the right-trivialized differential $\mathrm{dcay%
}_{\tilde{\mathbf{x}}}:\mathfrak{g}\rightarrow \mathfrak{g}$ in (\ref%
{diffCay}) is \cite{CelledoniOwren2003,HairerLubichWanner2006}%
\begin{equation}
\mathrm{dcay}_{\tilde{\mathbf{x}}}\left( \tilde{\mathbf{y}}\right) =2\left( 
\mathbf{I}-\tilde{\mathbf{x}}\right) ^{-1}\tilde{\mathbf{y}}\left( \mathbf{I}%
+\tilde{\mathbf{x}}\right) ^{-1}.  \label{dCay1}
\end{equation}%
Notice that the Cayley transform is occasional defined with $\tilde{\mathbf{x%
}}$ in (\ref{Cay1}) divided by 2. Then the factor 2 in (\ref{dCay1})
disappears. A motivation for doing so is that then (\ref{dCay1}) is the
identity mapping for $\mathbf{x}=\mathbf{0}$.

The inverse of the right trivialized differential is immediately found from (%
\ref{Cay2}) \cite{CelledoniOwren2003,HairerLubichWanner2006} 
\begin{align}
\mathrm{dcay}_{\tilde{\mathbf{x}}}^{-1}\left( \tilde{\mathbf{y}}\right) &= 
\frac{1}{2}\left( \mathbf{I}-\tilde{\mathbf{x}}\right) \tilde{\mathbf{y}}%
\left( \mathbf{I}+\tilde{\mathbf{x}}\right)  \notag \\
&=\frac{1}{2}\left( \tilde{\mathbf{y}}-\left[ \tilde{\mathbf{x}},\tilde{%
\mathbf{y}}\right] -\tilde{\mathbf{x}}\tilde{\mathbf{y}}\tilde{\mathbf{x}}%
\right)  \label{dCayInv1}
\end{align}%
with the matrix commutator $\left[ \tilde{\mathbf{x}},\tilde{\mathbf{y}}%
\right] =\tilde{\mathbf{x}}\tilde{\mathbf{y}}-\tilde{\mathbf{y}}\tilde{%
\mathbf{x}}$ as Lie bracket.

\subsection{Spatial Rotations -- $SO\left( 3\right) $%
%TCIMACRO{\TeXButton{secCaySO3}{\label{secCaySO3}}}%
%BeginExpansion
\label{secCaySO3}%
%EndExpansion
}

\subsubsection{Gibbs-Rodrigues parameters}

Rodrigues \cite{Rodrigues1840,AltmannBook1986} derived the rotation of a
vector as combination of two half-rotations about the rotation axis, which
leads to the corresponding rotation matrix. The latter can be constructed by
means of the Cayley transform on $so\left( 3\right) $, as a formalization of
Cayley's original result \cite{Cayley1843}. Invoking again the relation ${%
\tilde{\mathbf{x}}^{3}=}-\left\Vert \mathbf{x}\right\Vert ^{2}\tilde{\mathbf{%
x}}$, the Cayley map attains the well-known closed form \cite%
{BottemaRoth1979,McCarthyBook1990}%
\begin{align}
\mathrm{cay}\left( \tilde{\mathbf{x}}\right) & =\left( \mathbf{I}-\tilde{%
\mathbf{x}}\right) ^{-1}\left( \mathbf{I}+\tilde{\mathbf{x}}\right) =(%
\mathbf{I}+\tilde{\mathbf{x}}+\tilde{\mathbf{x}}^{2}+\tilde{\mathbf{x}}%
^{3}+\ldots )\left( \mathbf{I}+\tilde{\mathbf{x}}\right)  \notag \\
& =\mathbf{I}+\sigma (\tilde{\mathbf{x}}+\tilde{\mathbf{x}}^{2}).
\label{CaySO31}
\end{align}%
with%
\begin{equation}
\sigma :=\frac{2}{1+\left\Vert \mathbf{x}\right\Vert ^{2}}
\end{equation}%
Vector $\mathbf{x}\in {\mathbb{R}}^{3}$ is the \emph{Gibbs-Rodrigues vector}
while its components are the \emph{Rodrigues parameters}. The
Gibbs-Rodrigues vector is given in terms of the rotation angle $\varphi $
and the unit vector $\mathbf{n}$ along the rotation axis as $\mathbf{x}=\tan 
\frac{\varphi }{2}\tilde{\mathbf{n}}$. A singularity is encountered at $%
\varphi \rightarrow \pm \pi $ for which $\left\Vert \mathbf{x}\right\Vert
\rightarrow \infty $. Apparently rotations can be represented with $\varphi
\in \left( -\pi ,\pi \right) $ only, i.e. the Cayley map cannot be used to
describe full turn or even multiple turn rotations. This issue can be
tackled by using higher-order Cayley transformation \cite%
{SchaubTsiotrasJunkins1995,TsiotrasJunkinsSchaub1997,Tsiotras1998}, where
the singularities are shifted to multiple full turns. The conformal rotation
vector, also introduced by Milenkovic \cite{Milenkovic1982,Milenkovic2000}
and Wiener \cite[p. 69]{Wiener1962}, is a special case, which allows
describing double rotations. It should be recalled that Lie group
integration methods using the Cayley map solve (\ref{RecCay}), respectively
the right equation in (\ref{RecExpCay}), for the incremental rotation within
a time step with initial value $\mathbf{X}=\mathbf{0}$, \emph{assuming} that
no full rotation takes place within a time step.

The Gibbs-Rodrigues parameterization is related to the axis-angle
description via%
\begin{equation}
\exp (\varphi \tilde{\mathbf{n}})=\mathrm{cay}\left( \tan \frac{\varphi }{2}%
\tilde{\mathbf{n}}\right) =\mathbf{I}+\frac{2}{1+\tan ^{2}\frac{\varphi }{2}}%
\left( \tilde{\mathbf{n}}+\tilde{\mathbf{n}}^{2}\right) .
\end{equation}%
It can also be related to the description in terms of unit quaternions. If $%
Q=\left( q_{0},\mathbf{q}\right) \in Sp\left( 1\right) $ denotes a unit
quaternion, then the Gibbs-Rodrigues vector, describing the same rotation,
is obtained as $\mathbf{x}=\mathbf{q}/q_{0}$. Unit quaternions and
Gibbs-Rodrigues vectors are thus related via the gnomic projection of ${%
\mathbb{R}}^{3}$ onto the unit sphere $S^{3}\cong Sp\left( 1\right) $.

\subsubsection{Differential of the Cayley map}

Closed form expressions of the differential of the Cayley map on $SO\left(
3\right) $ were reported in \cite%
{BorriTrainelliBottasso2000,CelledoniOwren2003}. A proof of these relations
is given in the following.

\begin{lemma}
The right-trivialized differential of $\mathrm{cay}:so\left( 3\right)
\rightarrow so\left( 3\right) $ and its inverse possess the matrix form%
\begin{align}
\mathbf{dcay}_{\mathbf{x}}& =\sigma \left( \mathbf{I}+\tilde{\mathbf{x}}%
\right)  \label{dCaySO31} \\
\mathbf{dcay}_{\mathbf{x}}^{-1}& =\frac{1}{\sigma }\mathbf{I}+\frac{1}{2}(%
\tilde{\mathbf{x}}^{2}-\tilde{\mathbf{x}})  \label{dCayInvSO31} \\
& =\frac{1}{2\sigma }\left( \mathbf{I}+\mathrm{cay}\left( -\tilde{\mathbf{x}}%
\right) \right) =\frac{1}{2\sigma }(\mathbf{I}+\mathbf{R}^{T})
\label{dCayInvSO32}
\end{align}%
with $\mathbf{R}=\mathrm{cay}\left( \tilde{\mathbf{x}}\right) $. They
satisfy the properties $\mathbf{dcay}_{-\mathbf{x}}=\mathbf{dcay}_{\mathbf{x}%
}^{T}$ and $\mathbf{dcay}_{-\mathbf{x}}^{-1}=\mathbf{dcay}_{\mathbf{x}}^{-T}$%
.
\end{lemma}

\begin{proof}
The last term in (\ref{dCayInv1}) can be expressed as $\tilde{\mathbf{x}}%
\tilde{\mathbf{y}}\tilde{\mathbf{x}}=\tilde{\mathbf{x}}^{2}\tilde{\mathbf{y}}%
-\tilde{\mathbf{x}}\widetilde{\tilde{\mathbf{x}}\mathbf{y}}=\mathbf{x}^{T}%
\mathbf{x}\tilde{\mathbf{y}}-\tilde{\mathbf{x}}\widetilde{\tilde{\mathbf{x}}%
\mathbf{y}}-\left\Vert \mathbf{x}\right\Vert ^{2}\tilde{\mathbf{y}}=%
\widetilde{\tilde{\mathbf{x}}\mathbf{y}}\tilde{\mathbf{x}}-\tilde{\mathbf{x}}%
\widetilde{\tilde{\mathbf{x}}\mathbf{y}}-\left\Vert \mathbf{x}\right\Vert
^{2}\tilde{\mathbf{y}}=\left[ \tilde{\mathbf{x}},\left[ \tilde{\mathbf{x}},%
\tilde{\mathbf{y}}\right] \right] -\left\Vert \mathbf{x}\right\Vert ^{2}%
\tilde{\mathbf{y}}$, and hence $\mathrm{dcay}_{\tilde{\mathbf{x}}%
}^{-1}\left( \tilde{\mathbf{y}}\right) =\frac{1}{2}((1+\left\Vert \mathbf{x}%
\right\Vert ^{2})\tilde{\mathbf{y}}-\left[ \tilde{\mathbf{x}},\tilde{\mathbf{%
y}}\right] +\left[ \tilde{\mathbf{x}},\left[ \tilde{\mathbf{x}},\tilde{%
\mathbf{y}}\right] \right] )$. With $\left[ \tilde{\mathbf{x}},\tilde{%
\mathbf{y}}\right] =\widetilde{\tilde{\mathbf{x}}\mathbf{y}}$ and $\left[ 
\tilde{\mathbf{x}},\left[ \tilde{\mathbf{x}},\tilde{\mathbf{y}}\right] %
\right] =\widetilde{\tilde{\mathbf{x}}\tilde{\mathbf{x}}\mathbf{y}}$ follows
the matrix (\ref{dCayInvSO31}) so that $\mathbf{z}=\mathbf{dcay}_{\mathbf{x}%
}^{-1}\mathbf{y}$ when $\tilde{\mathbf{z}}=\mathrm{dcay}_{\tilde{\mathbf{x}}%
}^{-1}\left( \tilde{\mathbf{y}}\right) $, and noting (\ref{CaySO31}) yields (%
\ref{dCayInvSO32}). It is easy to show (\ref{dCaySO31}) by multiplication
with (\ref{dCayInvSO31}). The closed form expression (\ref{dCaySO31}) shows
that $\mathbf{dcay}_{-\mathbf{x}}=\mathbf{dcay}_{\mathbf{x}}^{T}$.
\end{proof}

The spatial and body-fixed angular velocity and the time derivative of the
Gibbs-Rodrigues vector (also called Cayley quasi-velocities \cite%
{SinclairHurtado2005}) are thus related as%
\begin{align}
%TCIMACRO{\TeXButton{w}{\bm{\omega }}}%
%BeginExpansion
\bm{\omega }%
%EndExpansion
^{\mathrm{s}}& =\mathbf{dcay}_{\mathbf{x}}\dot{\mathbf{x}},\ \ \dot{\mathbf{x%
}}=\mathbf{dcay}_{\mathbf{x}}^{-1}%
%TCIMACRO{\TeXButton{w}{\bm{\omega }}}%
%BeginExpansion
\bm{\omega }%
%EndExpansion
^{\mathrm{s}} \\
%TCIMACRO{\TeXButton{w}{\bm{\omega }}}%
%BeginExpansion
\bm{\omega }%
%EndExpansion
^{\mathrm{b}}& =\mathbf{dcay}_{\mathbf{x}}^{T}\dot{\mathbf{x}},\ \ \dot{%
\mathbf{x}}=\mathbf{dcay}_{\mathbf{x}}^{-T}%
%TCIMACRO{\TeXButton{w}{\bm{\omega }}}%
%BeginExpansion
\bm{\omega }%
%EndExpansion
^{\mathrm{b}}.  \notag
\end{align}

\subsection{Euclidean Motions -- $SE\left( 3\right) $%
%TCIMACRO{\TeXButton{secCaySE3}{\label{secCaySE3}}}%
%BeginExpansion
\label{secCaySE3}%
%EndExpansion
}

\subsubsection{Cayley map on $SE\left( 3\right) $}

The Cayley map for Euclidean motions is obtained with the general relation (%
\ref{Cay1}) applied to matrices (\ref{se3}). The explicit form of the power
of $\hat{\mathbf{X}}\in se\left( 3\right) $, shown in (\ref{expSE3}), yields%
\begin{align}
\mathrm{cay}(\hat{\mathbf{X}})& =\left( 
\begin{array}{cc}
\mathrm{cay}\left( \tilde{\mathbf{x}}\right) & \left( \mathbf{I}+\mathrm{cay}%
\left( \tilde{\mathbf{x}}\right) \right) \mathbf{y} \\ 
\mathbf{0} & 1%
\end{array}%
\right)  \label{CaySE31} \\
& =\left( 
\begin{array}{cc}
\mathrm{cay}\left( \tilde{\mathbf{x}}\right) & 2\sigma \,\mathbf{dcay}_{%
\tilde{\mathbf{x}}}^{-T}\mathbf{y} \\ 
\mathbf{0} & 1%
\end{array}%
\right)  \label{CaySE32}
\end{align}%
where (\ref{CaySE32}) is obtained from (\ref{CaySE31}) using (\ref%
{dCayInvSO32}). The form (\ref{CaySE31}) was presented in \cite%
{Selig-IFToMM2007}. There is no established name for the parameter vector $%
\mathbf{X}\in {\mathbb{R}}^{6}$ in $\mathrm{cay}(\hat{\mathbf{X}})$.

\subsubsection{Differential of the Cayley map}

The local reconstruction equations (\ref{RecCay}) in terms of the Cayley
map, with Rodrigues parameters as local coordinates, can be written in
vector form as ${\mathbf{V}}{^{\mathrm{s}}}=\mathbf{dcay}_{\mathbf{X}}\dot{%
\mathbf{X}}$ and ${\mathbf{V}}{^{\mathrm{b}}}=\mathbf{dcay}_{-\mathbf{X}}%
\dot{\mathbf{X}}$, respectively, where $\mathbf{dcay}_{\mathbf{X}}:{\mathbb{R%
}}^{6}\rightarrow {\mathbb{R}}^{6}$ is the coefficient matrix of this linear
relation.

\begin{lemma}
The coefficient matrix of the right-trivialized differential of $\mathrm{cay}%
:se\left( 3\right) \rightarrow SE\left( 3\right) $ and its inverse in vector
representation posses the closed form expressions%
\begin{align}
\mathbf{dcay}_{\mathbf{X}}& =\left( 
\begin{array}{cc}
\mathbf{dcay}_{\tilde{\mathbf{x}}} & \mathbf{0} \\ 
\sigma \tilde{\mathbf{y}}\left( \mathbf{I}+\tilde{\mathbf{x}}\right) & \ \ \
2\mathbf{I}+\sigma (\tilde{\mathbf{x}}+\tilde{\mathbf{x}}^{2})%
\end{array}%
\right) =\left( 
\begin{array}{cc}
\sigma \left( \mathbf{I}+\tilde{\mathbf{x}}\right) & \mathbf{0} \\ 
\sigma \tilde{\mathbf{y}}\left( \mathbf{I}+\tilde{\mathbf{x}}\right) & \ \ \
2\mathbf{I}+\sigma (\tilde{\mathbf{x}}+\tilde{\mathbf{x}}^{2})%
\end{array}%
\right)  \label{dCaySE3} \\
\mathbf{dcay}_{\mathbf{X}}^{-1}& =\left( 
\begin{array}{cc}
\mathbf{dcay}_{\tilde{\mathbf{x}}}^{-1} & \mathbf{0} \\ 
\frac{1}{2}\left( \tilde{\mathbf{x}}-\mathbf{I}\right) \tilde{\mathbf{y}} & 
\ \frac{1}{2}\left( \mathbf{I}-\tilde{\mathbf{x}}\right)%
\end{array}%
\right) =\frac{1}{2}\left( 
\begin{array}{cc}
\frac{2}{\sigma }\mathbf{I}+\tilde{\mathbf{x}}^{2}-\tilde{\mathbf{x}} & 
\mathbf{0} \\ 
\left( \tilde{\mathbf{x}}-\mathbf{I}\right) \tilde{\mathbf{y}} & \ \mathbf{I}%
-\tilde{\mathbf{x}}%
\end{array}%
\right) .  \label{dCayInvSE3}
\end{align}
\end{lemma}

\begin{proof}
Inserting matrices of the form (\ref{se3}) corresponding to $\mathbf{X}%
=\left( \mathbf{x},\mathbf{y}\right) $ and $\mathbf{U}=\left( \mathbf{u},%
\mathbf{v}\right) \in {\mathbb{R}}^{6}$ into (\ref{dCayInv1}) yields%
\begin{equation}
\mathrm{dcay}_{\hat{\mathbf{X}}}^{-1}(\hat{\mathbf{U}})=\left( 
\begin{array}{cc}
\mathrm{dcay}_{\tilde{\mathbf{x}}}^{-1}\left( \tilde{\mathbf{u}}\right) & \
\ \frac{1}{2}\left( \mathbf{I}-\tilde{\mathbf{x}}\right) \left( \mathbf{v}-%
\tilde{\mathbf{y}}\mathbf{u}\right) \\ 
\mathbf{0} & 0%
\end{array}%
\right) \in se\left( 3\right) .
\end{equation}%
Expressing $\hat{\mathbf{Z}}=\mathrm{dcay}_{\hat{\mathbf{X}}}^{-1}(\hat{%
\mathbf{U}})$ in vector form as $\mathbf{Z}=\mathbf{dcay}_{\hat{\mathbf{X}}%
}^{-1}\mathbf{U}$, along with (\ref{dCayInvSO31}), yields (\ref{dCayInvSE3}%
). The inverse of (\ref{dCayInvSE3}) is readily obtained as%
\begin{equation}
\mathbf{dcay}_{\hat{\mathbf{X}}}=\frac{2}{1+\left\Vert \mathbf{x}\right\Vert
^{2}}\left( 
\begin{array}{cc}
\left( \mathbf{I}+\tilde{\mathbf{x}}\right) & \mathbf{0} \\ 
\tilde{\mathbf{y}}\left( \mathbf{I}+\tilde{\mathbf{x}}\right) & \ \
(1+\left\Vert \mathbf{x}\right\Vert ^{2})\mathbf{I}+\tilde{\mathbf{x}}+%
\tilde{\mathbf{x}}^{2}%
\end{array}%
\right)  \label{dCaySE31}
\end{equation}%
and along with (\ref{dCaySO31}) yields (\ref{dCaySE3}).
\end{proof}

The closed form expressions for the right-trivialized differential and its
inverse were reported in \cite{KobilarovMarsden2011,KobilarovICRA2014}
without proof.

\subsection{Directional derivative of the right-trivialized differential}

If the Cayley map is to be used in generalized-$\alpha $ Lie group schemes,
the directional derivative of dcay is needed. This derived next.

\begin{lemma}
The directional derivative of $\mathbf{dcay}_{\mathbf{x}}:so\left( 3\right)
\rightarrow so\left( 3\right) $ and its inverse is given in closed form as%
\begin{align}
\left( \mathrm{D}_{\mathbf{x}}\mathbf{dcay}\right) 
%TCIMACRO{\TeXButton{-0.5ex}{\hspace{-0.5ex}}}%
%BeginExpansion
\hspace{-0.5ex}%
%EndExpansion
\left( \tilde{\mathbf{y}}\right) & =\sigma \tilde{\mathbf{y}}-\sigma ^{2}(%
\mathbf{x}^{T}\mathbf{y)}\left( \mathbf{I}+\tilde{\mathbf{x}}\right)
\label{DdcaySO3} \\
(\mathrm{D}_{\mathbf{x}}\mathbf{dcay}^{-1})%
%TCIMACRO{\TeXButton{-0.5ex}{\hspace{-0.5ex}}}%
%BeginExpansion
\hspace{-0.5ex}%
%EndExpansion
\left( \tilde{\mathbf{y}}\right) & =(\mathbf{x}^{T}\mathbf{y})\mathbf{I}+%
\frac{1}{2}\left( \tilde{\mathbf{x}}\tilde{\mathbf{y}}+\tilde{\mathbf{y}}%
\tilde{\mathbf{x}}-\tilde{\mathbf{y}}\right)  \label{DdcayInvSO3}
\end{align}%
The directional derivative of matrix $\mathbf{dcay}_{\mathbf{X}}$ of the
right-trivialized differential in (\ref{dCaySE3}) and its inverse admit the
explicit following forms, where $\mathbf{X}=\left( \mathbf{x},\mathbf{y}%
\right) $ and $\mathbf{U}=\left( \mathbf{u},\mathbf{v}\right) $,%
\begin{equation}
\left( \mathrm{D}_{\mathbf{X}}\mathbf{dcay}\right) 
%TCIMACRO{\TeXButton{-0.4ex}{\hspace{-0.4ex}}}%
%BeginExpansion
\hspace{-0.4ex}%
%EndExpansion
(\mathbf{U})=\left( 
%TCIMACRO{\TeXButton{-1.5ex}{\hspace{-1.5ex}}}%
%BeginExpansion
\hspace{-1.5ex}%
%EndExpansion
\begin{array}{cc}
\sigma \tilde{\mathbf{u}}-\sigma ^{2}(\mathbf{x}^{T}\mathbf{u})\left( 
\mathbf{I}+\widetilde{\mathbf{x}}\right) & \ \ \mathbf{0} \\ 
\sigma \left( \tilde{\mathbf{v}}+\tilde{\mathbf{v}}\tilde{\mathbf{x}}+\tilde{%
\mathbf{y}}\tilde{\mathbf{u}}\right) -\sigma ^{2}(\mathbf{x}^{T}\mathbf{u)}%
\left( \tilde{\mathbf{y}}+\tilde{\mathbf{y}}\widetilde{\mathbf{x}}\right) & 
\sigma \left( \tilde{\mathbf{u}}+\tilde{\mathbf{x}}\tilde{\mathbf{u}}+\tilde{%
\mathbf{u}}\tilde{\mathbf{x}}\right) -\sigma ^{2}(\mathbf{x}^{T}\mathbf{u})(%
\tilde{\mathbf{x}}+\widetilde{\mathbf{x}}^{2})%
\end{array}%
%TCIMACRO{\TeXButton{-1.5ex}{\hspace{-1.5ex}}}%
%BeginExpansion
\hspace{-1.5ex}%
%EndExpansion
\right) 
%TCIMACRO{\TeXButton{TeX field}{\vspace{-1ex}} }%
%BeginExpansion
\vspace{-1ex}
%EndExpansion
\label{DdcaySE3}
\end{equation}%
\begin{equation}
(\mathrm{D}_{\mathbf{X}}\mathbf{dcay}^{-1})%
%TCIMACRO{\TeXButton{-0.3ex}{\hspace{-0.3ex}}}%
%BeginExpansion
\hspace{-0.3ex}%
%EndExpansion
(\mathbf{U})=\frac{1}{2}\left( 
%TCIMACRO{\TeXButton{-1.2ex}{\hspace{-1.2ex}}}%
%BeginExpansion
\hspace{-1.2ex}%
%EndExpansion
\begin{array}{cc}
2(\mathbf{x}^{T}\mathbf{u})\mathbf{I}+\tilde{\mathbf{u}}\tilde{\mathbf{x}}+%
\tilde{\mathbf{x}}\tilde{\mathbf{u}}-\tilde{\mathbf{u}} & \ \ \mathbf{0} \\ 
\tilde{\mathbf{x}}\tilde{\mathbf{v}}+\tilde{\mathbf{u}}\tilde{\mathbf{y}}-%
\tilde{\mathbf{v}} & -\tilde{\mathbf{u}}%
\end{array}%
%TCIMACRO{\TeXButton{-1.2ex}{\hspace{-1.2ex}}}%
%BeginExpansion
\hspace{-1.2ex}%
%EndExpansion
\right) .  \label{DdcayInvSE3}
\end{equation}
\end{lemma}

\begin{proof}
The expressions (\ref{DdcaySO3}) and (\ref{DdcayInvSO3}) follow with (\ref%
{DD}) directly from (\ref{dCaySO31}) and (\ref{dCayInvSO31}), respectively.
Noting that the block entries of (\ref{dCaySE3}) and (\ref{dCayInvSE3})
depend on $\mathbf{X}$, the expressions (\ref{DdcaySE3}) and (\ref%
{DdcayInvSE3}) are readily found with (\ref{DD}).
\end{proof}

The expressions (\ref{DdcaySO3})-(\ref{DdcayInvSE3}) seem not to be present
in the current literature. They will be crucial for constructing iteration
matrices within numerical time stepping schemes such as the generalized-$%
\alpha $ method.

\subsection{Adjoint Representation -- The 'Configuration Tensor'}

Explicit forms of the Cayley map for the adjoint representation (\ref{AdSE3}%
) and its derivative were reported in \cite{BorriTrainelliBottasso2000},
were it is referred to as the configuration tensor. In the following, the
Cayley map for the adjoint representing is presented to emphasize that the
involved parameters are different from those of the Cayley map on $SE\left(
3\right) $. The Cayley map for the adjoint representation is formally
introduced as%
\begin{equation}
\mathrm{cay}(\mathbf{ad}_{\mathbf{X}})=\left( \mathbf{I}-\mathbf{ad}_{%
\mathbf{X}}\right) ^{-1}\left( \mathbf{I}+\mathbf{ad}_{\mathbf{X}}\right) .
\label{CayAd}
\end{equation}

\begin{lemma}
The Cayley map of the adjoint representation of $SE\left( 3\right) $ admits
the following closed form expressions%
\begin{equation}
\mathrm{cay}(\mathbf{ad}_{\mathbf{X}})=\left( 
\begin{array}{cc}
\mathrm{cay}\left( \tilde{\mathbf{x}}\right) & \mathbf{0} \\ 
\mathbf{A}\left( \mathbf{x},\mathbf{y}\right) & \mathrm{cay}\left( \tilde{%
\mathbf{x}}\right)%
\end{array}%
\right)  \label{CayAd2}
\end{equation}%
with 
\begin{align}
\mathbf{A}\left( \tilde{\mathbf{x}},\tilde{\mathbf{y}}\right) & =\frac{1}{2}%
\left( \mathbf{I}+\mathbf{R}\right) \tilde{\mathbf{y}}\left( \mathbf{I}+%
\mathbf{R}\right)  \label{CayAd5} \\
& =2\left( \mathbf{I}-\tilde{\mathbf{x}}\right) ^{-1}\tilde{\mathbf{y}}%
\left( \mathbf{I}-\tilde{\mathbf{x}}\right) ^{-1}  \label{CayAd6} \\
& =\left( \mathrm{D}_{\tilde{\mathbf{x}}}\mathrm{cay}\right) \left( \tilde{%
\mathbf{y}}\right)  \label{CayAd7} \\
& =\mathrm{dcay}_{\tilde{\mathbf{x}}}\left( \tilde{\mathbf{y}}\right) 
\mathbf{R}  \label{CayAd8} \\
& =\tilde{\mathbf{y}}~\mathbf{dcay}_{\mathbf{x}}+\sigma \tilde{\mathbf{x}}%
\tilde{\mathbf{y}}\mathbf{R}  \label{CayAd9}
\end{align}%
with $\mathbf{R}:=\mathrm{cay}\left( \tilde{\mathbf{x}}\right) $.
\end{lemma}

\begin{proof}
The first term in the product (\ref{CayAd}) is readily obtained with $%
\mathbf{ad}_{\mathbf{X}}$ in (\ref{adSE3}) as%
\begin{align}
\left( \mathbf{I}-\mathbf{ad}_{\mathbf{X}}\right) ^{-1}& =\mathbf{I}+\mathbf{%
ad}_{\mathbf{X}}+\mathbf{ad}_{\mathbf{X}}^{2}+\mathbf{ad}_{\mathbf{X}%
}^{3}+\ldots  \notag \\
& =\sum_{i=0}^{\infty }\left( 
\begin{array}{cc}
\tilde{\mathbf{x}}^{i} & \mathbf{0} \\ 
\mathbf{P}_{i}\left( \tilde{\mathbf{x}},\tilde{\mathbf{y}}\right) & \tilde{%
\mathbf{x}}^{i}%
\end{array}%
\right)
\end{align}%
with $\mathbf{P}_{i}$ in (\ref{P}). The expression (\ref{Cay3}) shows that $%
\left( \mathrm{D}_{\tilde{\mathbf{x}}}\mathrm{cay}\right) \left( \tilde{%
\mathbf{y}}\right) =2\sum_{i=0}^{\infty }\mathbf{P}\left( \tilde{\mathbf{x}},%
\tilde{\mathbf{y}}\right) =2\sum_{i=0}^{\infty }(\mathrm{D}_{\tilde{\mathbf{x%
}}}\tilde{\mathbf{x}}^{i})%
%TCIMACRO{\TeXButton{-0.5ex}{\hspace{-0.5ex}}}%
%BeginExpansion
\hspace{-0.5ex}%
%EndExpansion
\left( \tilde{\mathbf{y}}\right) $, and hence%
\begin{equation}
\left( \mathbf{I}-\mathbf{ad}_{\mathbf{X}}\right) ^{-1}=\sum_{i=0}^{\infty
}\left( 
\begin{array}{cc}
\left( \mathbf{I}-\tilde{\mathbf{x}}\right) ^{-1} & \mathbf{0} \\ 
\frac{1}{2}\left( \mathrm{D}_{\tilde{\mathbf{x}}}\mathrm{cay}\right) 
%TCIMACRO{\TeXButton{-0.5ex}{\hspace{-0.5ex}}}%
%BeginExpansion
\hspace{-0.5ex}%
%EndExpansion
\left( \tilde{\mathbf{y}}\right) & \left( \mathbf{I}-\tilde{\mathbf{x}}%
\right) ^{-1}%
\end{array}%
\right) .
\end{equation}%
Postmultiplication with $\mathbf{I}+\mathbf{ad}_{\mathbf{X}}$, according to (%
\ref{CayAd}), yields (\ref{CayAd2}) with 
\begin{align}
\mathbf{A}\left( \tilde{\mathbf{x}},\tilde{\mathbf{y}}\right) & =\frac{1}{2}%
\left( \mathrm{D}_{\tilde{\mathbf{x}}}\mathrm{cay}\right) 
%TCIMACRO{\TeXButton{-0.5ex}{\hspace{-0.5ex}}}%
%BeginExpansion
\hspace{-0.5ex}%
%EndExpansion
\left( \tilde{\mathbf{y}}\right) \left( \mathbf{I}+\tilde{\mathbf{x}}\right)
+\left( \mathbf{I}-\tilde{\mathbf{x}}\right) ^{-1}\tilde{\mathbf{y}}  \notag
\\
& =\left( \mathbf{I}-\tilde{\mathbf{x}}\right) \tilde{\mathbf{y}}(\mathbf{I}%
+\left( \mathbf{I}-\tilde{\mathbf{x}}\right) ^{-1}\left( \mathbf{I}+\tilde{%
\mathbf{x}}\right) )  \label{CayAd3} \\
& =\left( \mathbf{I}-\tilde{\mathbf{x}}\right) \tilde{\mathbf{y}}\left( 
\mathbf{I}+\mathbf{R}\right)  \label{CayAd4}
\end{align}%
where (\ref{CayAd3}) is obtained invoking (\ref{dCay1}). Since $\mathbf{R}$
and $\tilde{\mathbf{x}}$ are related via the Cayley transform, they satisfy $%
\mathbf{I}+\mathbf{R}=2\left( \mathbf{I}-\tilde{\mathbf{x}}\right) ^{-1}$.
Inserting this, and its inverse, into (\ref{CayAd4}) yields (\ref{CayAd5})
and (\ref{CayAd6}), respectively. Relation (\ref{Cay2}) shows the identity
of (\ref{CayAd6}) and (\ref{CayAd7}). Inserting (\ref{Cay4}) and (\ref{dCay1}%
) into (\ref{CayAd6}) leads to (\ref{CayAd8}), which is merely restating the
definition $\mathrm{dcay}_{\tilde{\mathbf{x}}}\left( \tilde{\mathbf{y}}%
\right) \mathrm{cay}\left( \tilde{\mathbf{x}}\right) =\left( \mathrm{D}_{%
\tilde{\mathbf{x}}}\mathrm{cay}\right) \left( \tilde{\mathbf{y}}\right) $.

Now (\ref{CayAd8}) can be written, with (\ref{CaySO31}) and (\ref{dCaySO31}%
), as%
\begin{align}
\mathrm{dcay}_{\tilde{\mathbf{x}}}\left( \tilde{\mathbf{y}}\right) \mathrm{%
cay}\left( \tilde{\mathbf{x}}\right) & =\left( \mathbf{dcay}_{\mathbf{x}}%
\mathbf{y}\right) ^{\sim }\mathrm{cay}\left( \tilde{\mathbf{x}}\right) 
\notag \\
& =\sigma \left( \left( \mathbf{I}+\tilde{\mathbf{x}}\right) \mathbf{y}%
\right) ^{\sim }\left( \mathbf{I}+\sigma (\tilde{\mathbf{x}}+\tilde{\mathbf{x%
}}^{2})\right)  \notag \\
& =\sigma \left( \tilde{\mathbf{y}}+\tilde{\mathbf{x}}\tilde{\mathbf{y}}-%
\tilde{\mathbf{y}}\tilde{\mathbf{x}}\right) \left( \mathbf{I}+\sigma (\tilde{%
\mathbf{x}}+\tilde{\mathbf{x}}^{2})\right)  \notag \\
& =\sigma \left( \tilde{\mathbf{y}}+\tilde{\mathbf{x}}\tilde{\mathbf{y}}-%
\tilde{\mathbf{y}}\tilde{\mathbf{x}}\right) +\sigma ^{2}\left( (1+\left\Vert 
\mathbf{x}\right\Vert ^{2})\tilde{\mathbf{y}}\tilde{\mathbf{x}}+\tilde{%
\mathbf{x}}\tilde{\mathbf{y}}(\tilde{\mathbf{x}}+\tilde{\mathbf{x}}%
^{2})\right)  \notag \\
& =\sigma \left( \tilde{\mathbf{y}}\left( \mathbf{I}+\tilde{\mathbf{x}}%
\right) +\tilde{\mathbf{x}}\tilde{\mathbf{y}}\left( \mathbf{I}+\sigma (%
\tilde{\mathbf{x}}+\tilde{\mathbf{x}}^{2})\right) \right)  \notag \\
& =\tilde{\mathbf{y}}\,\mathbf{dcay}_{\mathbf{x}}+\sigma \tilde{\mathbf{x}}%
\tilde{\mathbf{y}}\mathrm{cay}\left( \tilde{\mathbf{x}}\right)
\end{align}%
which proves (\ref{CayAd9}).
\end{proof}

The expression (\ref{CayAd6}) was derived in \cite{Selig-IFToMM2007}, and (%
\ref{CayAd8}) was presented in \cite{BorriTrainelliBottasso2000}. Relation (%
\ref{CayAd5}) was reported without proof in \cite{BauchauChoi2003}, and
elements of the parameter vector $\mathbf{X}\in {\mathbb{R}}^{6}$ were
called \emph{Cayley-Gibbs-Rodrigues motion parameters}.

Equation (\ref{CayAd8}) along with (\ref{AdSE3}) reveals that the position
determined by the Cayley map for the adjoint representation determines the
position $\mathbf{r}=\mathbf{dcay}_{\mathbf{x}}\mathbf{y}$, which shows a
striking similarity to the exponential map (\ref{expSE3}) on $SE\left(
3\right) $. On the other hand, the position vector determined by (\ref%
{CaySE31}) and (\ref{CaySE32}) is $\mathbf{r}=\left( \mathbf{I}+\mathrm{cay}%
\left( \tilde{\mathbf{x}}\right) \right) \mathbf{y}=2\sigma \,\mathbf{dcay}_{%
\tilde{\mathbf{x}}}^{-T}\mathbf{y}$. This reveals that the parameters $%
\mathbf{X}\in {\mathbb{R}}^{6}$ describing $\mathbf{C}\in SE\left( 3\right) $
and those describing the 'configuration tensor' $\mathbf{Ad}_{\mathbf{C}}$
are different. Consequently, the right-trivialized differential of the
Cayley map on $SE\left( 3\right) $ and of the adjoint representation are
different. This is due to the use of non-canonical coordinates. The
differential of the exponential map is indeed the same for both
representations (\ref{C}) and (\ref{adSE3}) in terms of canonical
coordinates.

\section{Conclusion%
%TCIMACRO{\TeXButton{secConclusion}{\label{secConclusion}}}%
%BeginExpansion
\label{secConclusion}%
%EndExpansion
}

Closed form relations for the exponential and the Cayley map, their
right-trivialized differentials, and the directional derivative of these
differentials play a crucial part in most Lie group integration schemes.
This includes the Munthe-Kaas methods as well the Lie group generalized-$%
\alpha $ method. The latter has become an alternative to classical
integration methods for multibody systems, in particular for systems
comprising flexible bodies undergoing large deformations. The Lie group
generalized-$\alpha $ method was originally derived in terms of canonical
coordinates, and thus involved the derivatives of the exponential on $%
SE\left( 3\right) $. Evidently, the Cayley map is computationally more
efficient, and can equally be used as coordinate map within this method.
This paper presents a comprehensive summary of all relevant closed form
expressions. The presented explicit relations for the Cayley map of
Euclidean motions, as well as its differential and directional derivative,
will facilitate the development of computationally efficient Lie group
generalized-$\alpha $ integration schemes in terms non-canonical local
coordinates. Besides these original relations, a novel derivation of the
closed form expressions of the trivialized derivative of the exponential on $%
SE\left( 3\right) $ were given.

\bibliographystyle{IEEEtran}
\bibliography{ExpCaySE3}

\end{document}